\documentclass[12pt]{article}
\usepackage{amsfonts}
\usepackage{amsmath}
\usepackage{amssymb}
\usepackage{amsthm}
\usepackage{a4wide}
\usepackage{todonotes}
\usepackage[utf8]{inputenc}

\newcommand{\abGal}[1]{\operatorname{Gal}\left(\overline{#1}/#1 \right)}
\newcommand{\ordPlaces}{\Omega_K^{\text{(ord)}}}

\newcommand{\EPlaces}{\Omega_K'}

\newcommand{\EndKBar}{\operatorname{End}^0_{\overline{K}}(A)}
\newcommand{\Kconn}{K^{\text{conn}}}
\title{Computing the geometric endomorphism ring of a genus 2 Jacobian}

\usepackage{scalerel,xcolor}
\def\darrow{\mathrel{\ThisStyle{\ooalign{$\SavedStyle\rightarrow$\cr%
  \hfil\textcolor{white}{\rule{2\LMpt}{1\LMex}}\kern2\LMpt\hfil}}}}

\usepackage[hyperfootnotes=false]{hyperref}

\newtheorem{theorem}{Theorem}

\newtheorem{algorithm}[theorem]{Algorithm}

\newtheorem{corollary}[theorem]{Corollary}

\newtheorem{lemma}[theorem]{Lemma}

\newtheorem{proposition}[theorem]{Proposition}

\theoremstyle{definition}
\newtheorem{remark}[theorem]{Remark}
\newtheorem{definition}[theorem]{Definition}
\newtheorem{example}[theorem]{Example}
\numberwithin{theorem}{section}

\newcommand\blfootnote[1]{%
  \begingroup
  \renewcommand\thefootnote{}\footnote{#1}%
  \addtocounter{footnote}{-1}%
  \endgroup
}

\date{}
\author{Davide Lombardo\thanks{Università di Pisa -- \texttt{davide.lombardo@unipi.it}}}
\begin{document}

\maketitle
\abstract{We describe an algorithm, based on the properties of the characteristic polynomials of Frobenius, to compute $\operatorname{End}_{\overline{K}}(A)$ when $A$ is the Jacobian of a nice genus-2 curve over a number field $K$. We use this algorithm to confirm that the description of the structure of the geometric endomorphism ring of $\operatorname{Jac}(C)$ given in the LMFDB ($L$-functions and modular forms database) is correct for all the genus 2 curves $C$ currently listed in it. We also discuss the determination of the field of definition of the endomorphisms in some special cases.}

\section{Introduction}
\blfootnote{MSC Classification: 11F80, 11G10, 11Y99}
\blfootnote{Keywords: abelian surfaces, Jacobian, Galois representations, endomorphisms}
The arithmetic study of Jacobians of genus-2 curves defined over the rationals (or, more generally, over number fields) from a computational point of view is a topic that has received much attention in recent times, but much still remains to be done. 
Let $K$ be a number field and $J$ be an abelian surface defined over $K$, which we will usually think of as being the Jacobian of a ``nice'' (that is, smooth projective) $K$-curve; while one could consider more general abelian surfaces, this is by far the most common case from the computational point of view. 
In this paper we consider the problem of provably determining the ring of endomorphisms of such a Jacobian, both from a theoretical and algorithmic point of view.

A brief outline of this article is as follows. We start in section \ref{sect:Preliminaries} by introducing the main actors of the paper: Galois representations, Mumford-Tate groups, characteristic polynomials of Frobenius, and the moduli space of genus-2 curves. In section \ref{sect:DeterminingItInPrinciple} we describe a night-and-day algorithm that is guaranteed to output the structure of $\operatorname{End}_K(A)$ for any abelian variety $A$ over a number field $K$, provided that one knows an explicit embedding of $A$ into projective space. Since this algorithm is hopelessly slow in practice, in section \ref{sect:ReallyDetermineEnd} we outline a different -- and this time computationally viable -- technique to compute $\EndKBar$ in the case of $A$ being a genus 2 Jacobians. This algorithm requires one to perform many sub-tasks, such as proving that $A$ is geometrically irreducible when this is the case (section \ref{sect_Irreducibility}), finding maps from $C$ to elliptic curves when they exist (section \ref{sect_Reducibility}), and determining $\EndKBar$ under the assumption that it is a field (section \ref{sect_Field}). The techniques described in section \ref{sect:ReallyDetermineEnd} mostly allow us to prove \textit{upper bounds} on $\EndKBar$ (namely, to show that it is contained in a certain $\mathbb{Q}$-algebra); we consider the opposite problem -- that of showing that $A$ admits nontrivial endomorphisms -- in section \ref{sect:Certify}. We show how to use the geometry of the moduli space of genus-2 curves to prove that $A=\operatorname{Jac}(C)$ admits real and complex multiplication (§\ref{sect:RM} and \ref{sect:CM}), and how to use results about real multiplication to also certify the existence of quaternionic multiplication (§\ref{sect:QM}). In section \ref{sect:PotentialRM} we also comment on the Galois structure of $\EndKBar$ when this is a real quadratic field.
Finally, in section \ref{sect:NumericalResults} we report on some numerical findings: our algorithms have allowed us to compute the structure of $\operatorname{End}_{\overline{\mathbb{Q}}}(\operatorname{Jac}(C))$ for all the genus-2 curves over $\mathbb{Q}$ admitting an odd-degree model with small coefficients and for all the curves considered in \cite{MR3540958}. In all cases our findings are in agreement with the data recorded in the \cite{lmfdb}.

To conclude this introduction, we point out to the reader that the recent preprint \cite{2017arXiv170509248C} develops techniques for the efficient computation of lower bounds -- that is, for the certification of the existence of extra endomorphisms -- which nicely complement the present paper. The approach of \cite{2017arXiv170509248C} relies heavily on the complex uniformization of abelian varieties, and is similar (though significantly more sophisticated) to that of section \ref{sect_Reducibility}. Of course, in order to completely certify the computation of $\EndKBar$ one also needs to prove upper bounds; this is what the present paper focuses on, and our results are used in part in \cite[Section 7]{2017arXiv170509248C} for precisely the purpose of showing that $\EndKBar$ is no larger than the analytic computations suggest. Another difference between the present work and \cite{2017arXiv170509248C} is that we also show how Frobenius polynomials can be used, at least in principle, to certify lower bounds.

\begin{remark}
Algorithms to verify the claims made in this paper are available at \url{http://people.dm.unipi.it/lombardo/}.
\end{remark}

\section{Preliminaries}\label{sect:Preliminaries}
\subsection{Notation}
Let $K$ be a number field and $A$ be an abelian surface over $K$. We write $\operatorname{End}^0_K(A)$ for the semisimple algebra $\operatorname{End}_K(A) \otimes_{\mathbb{Z}} \mathbb{Q}$.
We let $\Omega_K$ be the set of places of $K$. For $v \in \Omega_K$ we write $\mathbb{F}_v$ for the residue field at $v$, and let $p_v$ (resp.~$q_v$) be its residual characteristic (resp.~the cardinality of $\mathbb{F}_v$). We let $\Omega_K^{(A)}$ be the set of places at which $A$ has good reduction, and for $v \in \Omega_K^{(A)}$ we write $A_v$ for the reduction of $A$ at $v$ (an abelian variety over $\mathbb{F}_v$).

The natural action of $\abGal{K}$ on the Tate module $T_\ell(A)$ (respectively on the group $A[\ell^n]$ of $\ell^n$-torsion points) gives rise to a representation $\rho_{\ell^\infty}: \abGal{K} \to \operatorname{Aut} T_\ell(A)$ (respectively $\rho_{\ell^n}: \abGal{K} \to \operatorname{Aut} A[\ell^n]$) whose image we denote by $G_{\ell^\infty}$ (resp.~$G_{\ell^n}$). We shall also make use of the rational Tate module $V_\ell(A):=T_\ell(A) \otimes \mathbb{Q}_\ell$.
The symbol $\mathcal{G}_\ell$ will denote the $\mathbb{Q}_\ell$-Zariski closure of $G_{\ell^\infty}$ inside $\operatorname{Aut} V_\ell(A)$, while $\mathcal{G}_\ell^0$ will denote the connected component of the identity of $\mathcal{G}_\ell$. For each place $v$ at which $A$ has good reduction, we have a well-defined Frobenius polynomial $f_v(x) \in \mathbb{Z}[x]$ which, when $A$ is explicitly given as the Jacobian of a genus-2 curve, can be computed by counting points modulo $v$. Since places of degree 1 have full density in the set of places of any number field, for this point-counting one can restrict to working over prime fields. Moreover, for $K=\mathbb{Q}$ there are fast methods to compute the characteristic polynomial of Frobenius for all primes up to a large bound $N$: see \cite{MR3240808}, which gives an algorithm of time complexity $O(N^{1+\varepsilon})$.
We shall also find it useful to employ the following notation:

\begin{definition}
Let $F$ be any field and let $f(x) \in F[x]$ be a monic polynomial. Write $f(x)=\prod_{i=1}^n (x-x_i)$ for some $x_i \in \overline{F}$. For every positive integer $m$, we denote by $f^{[m]}$ the polynomial
\[
f^{[m]}(x)=\prod_{i=1}^n \left(x-x_i^m\right).
\]
\end{definition}
The following lemma is an obvious consequence of the fundamental theorem of symmetric functions, and implies that computing $f^{[m]}(x)$ from $f(x)$ for some fixed $m$ is an algorithmically easy matter:
\begin{lemma}
Let $F$ be a field, $f(x) \in F[x]$ be a polynomial of degree $n$, and $m$ be a positive integer. There exist polynomials $g_0,\ldots,g_n \in \mathbb{Z}[x_0,\ldots,x_n]$, depending only on $m$ and $n$, such that if $f(x)=a_nx^n+\cdots+a_0$ then
$
f^{[m]}(x)=g_n(a_0,\ldots,a_n)x^n+\cdots+g_0(a_0,\ldots,a_n).
$
\end{lemma}
\begin{remark}\label{rmk_CharPoly}
Suppose $f(x)$ as in the previous lemma is obtained as the characteristic polynomial of a certain endomorphism $g \in \operatorname{GL}_{n}(F)$: then $f^{[m]}(x)$ is the characteristic polynomial of $g^m$ (indeed both polynomials share the same roots and the same leading term).
\end{remark}
\begin{remark}\label{rmk:ReducibilityTwistedPolys}
Notice that $(fg)^{[a]}=f^{[a]}g^{[a]}$ and $f^{[ab]}=(f^{[a]})^{[b]}$. From this it follows that if $a \mid b$ and $f^{[a]}$ is reducible, then $f^{[b]}$ is also reducible.
\end{remark}

\subsection{Mumford-Tate groups}\label{sect_MT}

Let $K$ be a number field and $A$ be an abelian variety of dimension $g$ over $K$. We fix once and for all an embedding $\sigma$ of $K$ in $\mathbb{C}$ and we denote by $A_\mathbb{C}$ the base-change of $A$ to $\mathbb{C}$ along $\sigma$. The $\mathbb{Q}$-vector space $V=H_1 \left(A_{\mathbb{C}}(\mathbb{C}),\mathbb{Q}\right)$ is naturally endowed with a Hodge structure of type $(-1,0)\oplus(0,-1)$, that is, a decomposition of $\mathbb{C}$-vector spaces $V \otimes \mathbb{C} \cong V^{-1,0} \oplus V^{0,-1}$ such that $\overline{V^{-1,0}}=V^{0,-1}$.
Let $\mu_\infty:\mathbb{G}_{m,\mathbb{C}} \to \operatorname{GL}_{V \otimes \mathbb{C}}$ be the unique cocharacter such that $z \in \mathbb{C}^*$ acts as multiplication by $z$ on $V^{-1,0}$ and trivially on $V^{0,-1}$. 
\begin{definition}
The \textbf{Mumford-Tate} group of $A$ is the $\mathbb{Q}$-Zariski closure of the image of $\mu_\infty$, that is to say the smallest $\mathbb{Q}$-algebraic subgroup $\operatorname{MT}(A)$ of $\operatorname{GL}_V$ such that $\mu_\infty$ factors through $\operatorname{MT}(A)_\mathbb{C}$. 
\end{definition}

It is not hard to show that $\operatorname{MT}(A)$ is connected and contains the torus of homotheties in $\operatorname{GL}_V$.
As a consequence of the existence of polarizations on abelian varieties, it is also known that the Mumford-Tate group of an abelian variety of dimension $g$ is contained in the general symplectic group $\operatorname{GSp}_{2g,\mathbb{Q}}$. The precise symplectic form that is preserved (up to similitude) depends on the choice of a polarization, but (for our case of interest $g=2$) we can choose our basis of $V \cong \mathbb{Q}^{4}$ in such a way that the bilinear form preserved up to similitude is the one associated with the matrix
$
\begin{pmatrix}
0 & 1 & 0 & 0 \\
-1 & 0 & 0 &0 \\
0 & 0 & 0 & 1 \\
0 & 0 & -1 & 0
\end{pmatrix}
$. 
Recall now that we denote by $\mathcal{G}_\ell$ the Zariski closure of $G_{\ell^\infty}$, the image of Galois, in $\operatorname{Aut}\left(T_\ell(A) \otimes \mathbb{Q}_\ell\right)$. The interest of the Mumford-Tate group in our context is due to the celebrated Mumford-Tate conjecture, which predicts that (if one identifies $H_1(A_\mathbb{C}(\mathbb{C}),\mathbb{Q}) \otimes \mathbb{Q}_\ell$ with $V_\ell(A)$ via the comparison isomorphism of étale cohomology) the identity component $\mathcal{G}_\ell^0$ of $\mathcal{G}_\ell$ should coincide with $\operatorname{MT}(A) \times_\mathbb{Q} \mathbb{Q}_\ell$.
When $A$ is an abelian surface the Mumford-Tate conjecture is known to hold, so $\mathcal{G}_\ell^0$ is determined by $\operatorname{MT}(A)$ thanks to the canonical isomorphism $\mathcal{G}_\ell^0 = \operatorname{MT}(A) \times_{\mathbb{Q}} \mathbb{Q}_\ell$. In order to fully describe $\mathcal{G}_\ell$ we also need to understand the group of connected components $\mathcal{G}_\ell/\mathcal{G}_\ell^0$, which is known to be independent of $\ell$ thanks to the following theorem of Serre (see \cite[Proposition 6.14]{MR1150604} for a published account of the proof):
\begin{theorem}\label{thm:K'}
The kernel of the canonical morphism $\varepsilon_\ell : \abGal{K} \to \mathcal{G}_\ell/\mathcal{G}_\ell^0$ is independent of $\ell$. In particular, for all primes $p$ the fixed field $K'$ of $\ker \varepsilon_p$ is the unique minimal (and automatically normal) field extension of $K$ with the property that the Zariski closure of $\rho_{\ell^\infty}(\abGal{K'})$ is connected for all primes $\ell$. 
\end{theorem}

\begin{definition}
Let $A/K$ be an abelian variety. We denote by $\Kconn$ the field $K'$ whose existence is guaranteed by Theorem \ref{thm:K'}. It is in particular a \textit{normal} extension of $K$.
\end{definition}

Depending on the endomorphism algebra of $A$, we see from \cite{MR2982436} that the following are the only possibilities for $\operatorname{MT}(A)$ and $\mathcal{G}_\ell/\mathcal{G}_\ell^0$:
\begin{enumerate}
\item $\operatorname{End}^0_{\overline{K}}(A)=\mathbb{Q}$. The Mumford-Tate group of $A$ is $\operatorname{GSp}_{4,\mathbb{Q}}$, and the group of connected components $\mathcal{G}_\ell/\mathcal{G}_\ell^0$ is trivial.
\item $\operatorname{End}^0_{\overline{K}}(A)=F$, a real quadratic field. The Mumford-Tate group of $A$ is the quasi-direct product $\mathbb{G}_m \cdot \operatorname{Res}_{F/\mathbb{Q}}\operatorname{Sp}_{2,F}$, and the group of connected components $\mathcal{G}_\ell/\mathcal{G}_\ell^0$ is either trivial or of order 2.
\item $\operatorname{End}^0_{\overline{K}}(A)=F$, a CM field of degree 4. The Mumford-Tate group of $A$ is a certain (explicit) rank-3 subtorus of $\operatorname{Res}_{F/\mathbb{Q}}(\mathbb{G}_{m,F})$, and the group of connected components $\mathcal{G}_\ell/\mathcal{G}_\ell^0$ is isomorphic to one of the following 4 groups:
$\{1\}$, $C_2$, $C_4$, $C_2 \times C_2$.

\item $\operatorname{End}^0_{\overline{K}}(A)=Q$, a non-split quaternion algebra over $\mathbb{Q}$. The Mumford-Tate group of $A$ is $Q^\times$ (the groups of units of $Q$, considered as an algebraic group over $\mathbb{Q}$) and in particular has rank 2. The group of connected components $\mathcal{G}_\ell/\mathcal{G}_\ell^0$ is isomorphic to one of the following 9 groups:
$
\{1\}, C_2, C_3, C_4, C_6, D_2, D_3, D_4, D_6.
$

\item $\operatorname{End}^0_{\overline{K}}(A)=F_1 \oplus F_2$, where each $F_i$ is either $\mathbb{Q}$ or an imaginary quadratic field. In this case $A$ is geometrically isogenous to the product of two non-isogenous elliptic curves $E_1$ and $E_2$ that furthermore satisfy $\operatorname{End}^0_{\overline{K}}(E_i)=F_i$. The Mumford-Tate group has rank 3, and it is isomorphic to $\mathbb{G}_m \cdot \left( M_1 \times M_2 \right)$, where
\[
M_i = \begin{cases} \operatorname{SL}_{2,\mathbb{Q}}, \text{ if }F_i=\mathbb{Q} \\ \left\{ x \in \operatorname{Res}_{F_i/\mathbb{Q}}(\mathbb{G}_m) \bigm\vert x\overline{x}=1 \right\}, \text{ if }F_i\text{ is imaginary quadratic} \end{cases}
\]
The group of connected components $\mathcal{G}_\ell/\mathcal{G}_\ell^0$ is isomorphic to one of the following 4 groups:
$
\{1\}, C_2, C_4, C_2 \times C_2.
$
\item $\operatorname{End}^0_{\overline{K}}(A)=M_2(F)$, where $F$ is either $\mathbb{Q}$ or an imaginary quadratic field. In this case $A$ is geometrically isogenous to the square of an elliptic curve $E$ such that $\operatorname{End}^0_{\overline{K}}(E)=F$. The Mumford-Tate group is $\mathbb{G}_m \cdot \{ (x,x) \bigm\vert x \in M \}$, where as above $M$ is either $\operatorname{SL}_{2,\mathbb{Q}}$, if $F=\mathbb{Q}$, and $\left\{ x \in \operatorname{Res}_{F/\mathbb{Q}}(\mathbb{G}_m) \bigm\vert x\overline{x}=1 \right\}$, if $F$ is imaginary quadratic.
The group of connected components $\mathcal{G}_\ell/\mathcal{G}_\ell^0$ is isomorphic to a subgroup of either $S_4 \times C_2$ or $D_6 \times C_2$.
\end{enumerate}
Moreover, the minimal field of definition of the endomorphisms coincides with the field $K'$ of theorem \ref{thm:K'}.
Notice that in all cases the exponent of the group $\mathcal{G}_\ell/\mathcal{G}_\ell^0$ divides 12; in particular,
\begin{lemma}\label{lemma:ExpDivides12}
Let $A/K$ be an abelian surface and let $K'$ be the minimal extension of $K$ such that the Galois representations associated with $A/K'$ have connected image. 
\begin{enumerate}
\item All the endomorphisms of $A$ are defined over $K'$.
\item Let $w$ be a place of $K'$ lying above a place $v$ of $K$. The degree of the extension $\mathbb{F}_w/\mathbb{F}_v$ belongs to the set $\{1,2,3,4,6,12\}$.
\end{enumerate}
\end{lemma}

\subsection{Ordinary reduction}

We write $\ordPlaces$ for the set of the places of $K$ at which $A$ has good ordinary reduction. Recall that an abelian variety $A$ defined over a finite field $\mathbb{F}$ of characteristic $p$ is \textit{ordinary} if $\#A_{\overline{\mathbb{F}}}[p]=p^{\dim A}$ (see \cite[§3]{MR0314847} for more details). The following facts are well-known (for the first statement see for example \cite[Proposition 3.1]{MR1628150}; for the second, \cite[Theorem 8 and §3]{MR0314847}):

\begin{lemma}\label{lemma_ordImpliesEndDefBase}
Let $v$ be a place of characteristic $p$ at which $A$ has good reduction. Then $A_v$ is ordinary if and only if the characteristic polynomial of Frobenius $f_v(x)=x^4+a x^3+bx^2+apx+p^2$ satisfies $b \not \equiv 0 \pmod{p}$.
For all $v \in \ordPlaces$ such that $A_v$ is absolutely simple and for all positive integers $N$ we have 
$
\mathbb{Q}(\pi_v) = \mathbb{Q}(\pi_v^N) \cong  \operatorname{End}^0_{\mathbb{F}_v}(A_v) =\operatorname{End}^0_{\overline{\mathbb{F}_v}}(A_v),
$
where $\pi_v$ is a root of the (irreducible) polynomial $f_v(x)$ or, equivalently, the Frobenius automorphism of $A_v$.
\end{lemma}

We shall make use of the fact that ordinarity is clearly a \textit{geometric} property, hence independent of the field of definition. We record in the following corollary two consequences of this remark:
\begin{corollary}\label{cor_ordImpliesEndDefBase}
The following hold:
\begin{enumerate}
\item let $A$ be an abelian variety over a number field $K$ and let $K'$ be a finite extension of $K$. Suppose that $w$ is a place of $K'$ such that $A_w$ is absolutely simple and ordinary. Then $A_v$ is also absolutely simple and ordinary, where $v$ is the place of $K$ induced by $w$;
\item in the situation above, one has $\operatorname{End}^0_{\mathbb{F}_v}(A_v)=\operatorname{End}^0_{\mathbb{F}_w}(A_w)$.
\end{enumerate}
\end{corollary}
\begin{proof}
Let $\mathbb{F}_2/\mathbb{F}_1$ be an extension of finite fields and let $A'$ be an abelian variety defined over $\mathbb{F}_1$. If $A'_{\mathbb{F}_2}$ is ordinary, then $A'/\mathbb{F}_1$ is ordinary, because the definition of ordinarity only depends on $A'_{\overline{\mathbb{F}_1}} \cong A'_{\overline{\mathbb{F}_2}}$. (1) clearly follows; as for (2), using the fact that $A_v$ is ordinary (by (1)) we know from lemma \ref{lemma_ordImpliesEndDefBase} that $\operatorname{End}^0_{\mathbb{F}_v}(A_v)=\operatorname{End}^0_{\overline{\mathbb{F}_v}}(A_v)=\operatorname{End}^0_{\overline{\mathbb{F}_w}}(A_w)=\operatorname{End}^0_{\mathbb{F}_w}(A_w)$.
\end{proof}

 We shall also need the following fact, proved in \cite[Theorem 3, Corollary 2]{2015arXiv150604784S} as a consequence of the results of \cite{MR2982436}:
\begin{proposition}\label{prop:DensityOrdinaryReduction} Let $A/K$ be an abelian surface. The set of places of $K$ where $A$ has ordinary reduction admits natural density, and this density is either $1$, $1/2$ or $1/4$. Moreover, the set of places of $K$ where $A$ has ordinary reduction has density 1 in the following cases: $K \supseteq K^{\text{conn}}$; $\EndKBar$ is a real quadratic field; $\EndKBar$ is a nonsplit quaternion algebra.
\end{proposition}

\subsection{A theorem of Zarhin}

We now recall a result due to Zarhin, which will enable us to understand the $\ell$-adic completion of the endomorphism ring of the reductions $A_v$, at least for $v$ in a positive-density set of places of $K$.

\begin{theorem}[Zarhin \cite{MR3660686}]
\label{thm_Zarhin}
Let $A$ be an abelian variety of positive dimension over a number field $K$.
 Suppose that the groups $\mathcal{G}_{\ell}$ attached to $A$ are connected. Let $\mathbb{P}$ be a finite nonempty set of primes and suppose that for each $\ell \in \mathbb{P}$ we are given an element
$f_{\ell} \in \mathcal{G}_{\ell}(\mathbb{Q}_{\ell})\subset  \operatorname{Aut}_{\mathbb{Q}_{\ell}}(V_{\ell}(A))$
such that the characteristic polynomial 
$P_{f_{\ell}}(t)=\det (t \cdot \operatorname{Id} -f_{\ell} \bigm\vert V_{\ell}(A))\in \mathbb{Q}_{\ell}[t]$
has no multiple roots. Let
$\mathfrak{z}(f_{\ell})_0\subset \operatorname{End}_{\mathbb{Z}_{\ell}}(T_{\ell}(A))$
be the centralizer of $f_{\ell}$ in
$\operatorname{End}_{\mathbb{Z}_{\ell}}(T_{\ell}(A)) \subset \operatorname{End}_{\mathbb{Q}_{\ell}}(V_{\ell}(A)).$
Then the set of nonarchimedean places $v$ of $K$ such that 
the residual characteristic $p_v$ does not belong to $\mathbb{P}$, the abelian variety 
$A$ has good reduction at $v$ and
$$\operatorname{End}(A_v)\otimes \mathbb{Z}_{\ell} \cong \mathfrak{z}(f_{\ell})_0 \quad \forall \ell\in \mathbb{P}$$
has positive density. In addition, for all such $v$ the ring $\operatorname{End}(A_v)$ is commutative.
\end{theorem}

We will use this result in section \ref{sect_Field} to show that if $\operatorname{End}^0_{\overline{\mathbb{Q}}}(A)$ is a number field $E$, then one can determine $\operatorname{disc}(E)$ by looking only at the fields generated over $\mathbb{Q}$ by the roots of the characteristic polynomials of Frobenius. More precisely, theorem \ref{thm_Zarhin} will enable us to establish a relation between the $p$-parts of $\operatorname{disc}(E)$ and of $\operatorname{disc}(\operatorname{End}(A_v))$ which holds at least for a positive-proportion set of places $v$. Doing this for all primes $p$ will be enough to determine $\operatorname{disc}(E)$.

\subsection{The moduli space of genus-2 curves}\label{sect:ModuliSpace}
As it is well-known, the moduli space $\mathcal{M}_2$ of genus-2 curves over $\overline{\mathbb{Q}}$ is birationally equivalent to affine 3-space, the birational isomorphism being given by the so-called absolute Igusa invariants of the curve (see \cite{MR0114819}). Via the Torelli morphism, $\mathcal{M}_2$ can also be identified with the moduli space $\mathcal{A}_2$ of principally polarized abelian surfaces which are not products of elliptic curves. 
We shall interchangeably work in $\mathcal{A}_2$ or in its compactification $\mathcal{A}_2^* \cong \mathbb{P}_{\overline{\mathbb{Q}}}^3\left(I_2, I_4, I_6, I_{10}\right)$, where $I_2,\ldots,I_{10}$ are also called Igusa invariants (and where $I_{2k}$ has weight $2k$). One has
\[
\mathcal{A}_2 = \mathcal{A}_2^* \setminus H_1 \cong \mathbb{P}^3(I_2,I_4,I_6,I_{10}) \setminus \{I_{10}=0\},
\]
with the divisor at infinity $I_{10}=0$ corresponding to the locus of principally polarized abelian surfaces that are products of elliptic curves (with the product polarization).

In order to better describe the integral structure of the endomorphism rings we are interested in, we introduce the following general notion of optimal embedding:
\begin{definition}\label{def:QuadQuat}
An embedding $\rho: R \hookrightarrow S$ of rings is said to be \textbf{optimal} if the equality
\[
\{s \in S : \exists n \in \mathbb{Z}, ns \in \rho(R) \}=\rho(R)
\]
holds, or equivalently if $S/\rho(R)$ is torsion-free as an additive group. If $A$ is an abelian variety defined over a field $K$ and $R$ is a ring, we say that there is an \textbf{optimal action} of $R$ on $A$ defined over $K$ if there is an optimal embedding $R \hookrightarrow \operatorname{End}_K(A)$. When the field $K$ is not specified we mean that there exists an optimal embedding $R \hookrightarrow \operatorname{End}_{\overline{K}}(A)$.
\end{definition}

We shall also need the notion of discriminant and of quadratic  and quaternionic rings:
\begin{definition}\phantomsection\label{def:Discriminant}
\begin{itemize}
\item
A positive integer $D$ is a \textit{discriminant} if it is of the form $n^2d$ with $n,d$ positive integers, $d \equiv 0,1 \pmod 4$. When $n=1$, such a $D$ is called a \textit{fundamental discriminant}.
\item
For every positive discriminant $D$, there exists a unique (up to isomorphisms) ring $\mathcal{O}_D$ with the following property. There is an isomorphism of $\mathbb{Z}$-modules $\mathcal{O}_D \cong \mathbb{Z} \cdot 1 \oplus \mathbb{Z} \cdot \omega$ (for some $\omega \in \mathcal{O}_D$) such that, writing $\omega^2=a \omega+b$, one has $4D=a^2+4b$. When the ring $\mathcal{O}_D$ is an integral domain (i.e. if and only if $D$ is not a square), it is an order in the real quadratic field $\mathbb{Q}(\sqrt{D})$. Any such ring $\mathcal{O}_D$ is called a \textit{quadratic ring}.
\item A \textit{quaternionic ring} is an order in a quaternion algebra over $\mathbb{Q}$; as an abelian group under addition, any such ring is isomorphic to $\mathbb{Z}^4$.
\end{itemize}
\end{definition}

Jacobians with nontrivial endomorphisms are parametrized by proper subvarieties of $\mathcal{A}_2^*$ as follows (see for example \cite{MR1329525}); notice that $\mathcal{A}_2^*$ is the moduli space of Jacobians over $\overline{\mathbb{Q}}$, so in particular in each of the following cases the endomorphisms considered are defined over $\overline{\mathbb{Q}}$.
\begin{enumerate}
\item For every positive discriminant $D$ there is a (``Humbert'') hypersurface $H_{D}$ in $\mathcal{A}_2^*$ that parametrizes curves $C$ whose Jacobians admit an optimal action of $\mathcal{O}_D$. 
Humbert surfaces are irreducible and connected.

Points on $H_{n^2}$ parametrize curves whose Jacobian $J$ admits an $(n,n)$-isogeny to a product $E_1 \times E_2$ of two elliptic curves, and no isogeny of lower degree exists between $J$ and a product of two elliptic curves \cite{MR1285957}. Notice in particular that this is consistent with denoting by $H_1$ the locus of abelian surfaces that are products of elliptic curves.
\item For each quaternionic ring $R$ there are irreducible curves $S_{R,1}, \ldots, S_{R,k}$ contained in $\mathcal{A}_2^*$ that parametrizes curves whose Jacobians admit an optimal action of $R$. Different Shimura curves corresponding to the same $R$ parametrize different embeddings of $R$ in $\operatorname{M}_4(\mathbb{Z})$ (here we view $\operatorname{M}_4(\mathbb{Z})$ as acting on the rank-4 lattice that defines $J$ in $\mathbb{C}^2$).

We shall refer to each $S_{R,i}$ as a ``Shimura curve". We shall describe these quaternionic rings and Shimura curves in greater detail below in section \ref{sect:QuatRings}.
\item Finally, curves whose Jacobians admit complex multiplication correspond to isolated points in moduli space. Notice that some of these points lie on Shimura curves, and correspond to Jacobians that are isogenous to squares of elliptic curves with CM.
\end{enumerate}

From this description we have in particular:
\begin{proposition}
Let $J=\operatorname{Jac}(C)$ be the Jacobian of a curve whose corresponding point in moduli space is $x_J$. Then $J$ is a geometrically simple abelian variety if and only if $x_J \not \in \bigcup_{n \geq 1} H_{n^2}$.
\end{proposition}

\subsection{Quaternionic rings and Shimura curves}\label{sect:QuatRings}
Let $R$ be an order in a quaternion algebra $Q$ over $\mathbb{Q}$, and recall that any such $Q$ admits a canonical anti-involution, which we denote by $x \mapsto x^\dagger$. We define the reduced norm and trace of $x \in Q$ by $\operatorname{Trd}(x)=x+x^\dag$ and $\operatorname{Nrd}(x)=x x^\dagger$. Notice that $R$ contains a copy of $\mathbb{Z}$ which is given precisely by the set of fixed points of the canonical anti-involution restricted to $R$; it follows that (upon restriction to $R$) the reduced norm and trace give rise to maps $R \to \mathbb{Z}$.
The discriminant $d(x_1,x_2,x_3,x_4)$ of a 4-tuple $x_1,\ldots,x_4$ of elements of $R$ is an integer that satisfies
\[
d(x_1,x_2,x_3,x_4)^2 = - \det(\operatorname{Trd}(x_ix_j));
\]
conventions on the sign vary (as the discriminant is really an ideal and not an integer), but following \cite{MR1707758} we shall always choose the positive sign for $d(x_1,\ldots,x_4)$. The discriminant of $R$ is by definition the discriminant of a $\mathbb{Z}$-basis of $R$, and the discriminant of $Q$ is the discriminant of a maximal order in $Q$. One can check that these definitions depend neither on the choice of the basis in $R$ nor on that of the maximal order in $Q$.

For $x \in R$ we define the discriminant $\Delta(x):=\operatorname{Trd}(x)^2-4\operatorname{Nrd}(x)$, and we let 
\[
\Delta(x,y) := \frac{1}{2} \left( \Delta(x+y) - \Delta(x) - \Delta(y) \right);
\]
be the corresponding bilinear form.
Now let $R$ be a quaternionic ring which arises as the endomorphism ring of a genus-2 Jacobian. Any such ring admits a \textit{polarization} $\mu$, an element $\mu \in R$ such that $\mu^2 \in \mathbb{Z}, \mu^2<0$, and the map $x \mapsto \mu^{-1} x^\dag \mu$ is a positive involution on $R \otimes \mathbb{Q}$ (see \cite[§4.2.1]{Gruenewald} for more details). We call the pair $(R, \mu)$ a polarized quaternionic ring.

With each polarized quaternionic ring we can associate canonically an equivalence class of binary integral quadratic forms through the following procedure. There exist elements $\alpha, \beta \in R$ such that $\alpha\beta-\beta\alpha=-\mu$ and $R \cong \mathbb{Z} \oplus \mathbb{Z} \alpha \oplus \mathbb{Z} \beta \oplus \mathbb{Z} \alpha \beta$ as additive groups (see \cite[Proposition 4.4.1]{RotgerThesis}). The equivalence class of quadratic forms attached to $R$ is by definition the one containing the quadratic form $\Delta$, whose associated matrix is
\[
M_R =\begin{pmatrix}
\Delta(\alpha, \alpha) & \Delta(\alpha,\beta) \\ \Delta(\alpha,\beta) & \Delta(\beta,\beta)
\end{pmatrix}.
\]
Notice that $M_R(m,n)=\Delta(m\alpha+n\beta)=\operatorname{Trd}(m\alpha+n\beta)^2-4\operatorname{Nrd}(m\alpha+n\beta)$ is always a discriminant in the sense of definition \ref{def:Discriminant}.

The matrix $M_R$ encodes a number of useful informations about the ring $R$ itself:
\begin{theorem}{(\cite[Theorems 7 and 10 and Corollary 9]{MR1707758})}\label{thm_RMImpliesEverything} Let $R, M_R$ be as above, and let $S_R$ be a Shimura curve corresponding to $R$.
\begin{enumerate}
\item If $M_{R'}$ is another binary integral quadratic form arising from a quaternionic ring $R'$, then $M_R$ and $M_{R'}$ are $\operatorname{GL}_2(\mathbb{Z})$-equivalent if and only if $R$ and $R'$ are isomorphic as $\mathbb{Z}$-algebras.
\item $M_R$ is positive-definite.
\item $\det M_R = 4\operatorname{disc}(R)$.
\item Let $D$ be a positive discriminant. There is an optimal embedding of $\mathcal{O}_D$ in $R$ if and only if $M_R$ represents $D$ \textit{primitively}, that is, if and only if there exist integers $m, n$ with $(m,n)=1$ such that $M_R(m,n)=D$.
\item Suppose that $S_R$ is contained in the intersection $H_{D_1} \cap H_{D_2}$ of two distinct Humbert surfaces. Then for a suitable choice of basis of $R$ we have $M_R = \begin{pmatrix}
D_1 & k \\ k & D_2
\end{pmatrix}$ for some integer $k$.
\end{enumerate}
\end{theorem}

Notice that if $J$ is a genus-2 Jacobian, and if $R$ is the endomorphism ring of $J$, then there is an optimal embedding of $\mathcal{O}_D$ in $R=\operatorname{End}_{\overline{K}}(J)$ if and only if the point in moduli space corresponding to $J$ lies on the Humbert surface $H_D$.

Concerning the moduli interpretation of Shimura curves, we also remark that a point lying on the intersection of two Humbert surfaces $H_{D_1} \cap H_{D_2}$ with $D_1 \neq D_2$ corresponds either to a \textit{simple} abelian surface with quaternionic multiplication by an (automatically indefinite) quaternion algebra over $\mathbb{Q}$, or to the square of an elliptic curve (\cite[Proposition 2.15]{Gruenewald}). This is in particular true for points lying on Shimura curves.

\subsection{Some effective algorithms}\label{sect:BasicAlgorithms}
For the purposes of this section, we say that we are given an abelian variety $A$ over a number field $K$ if we know the set $S=\Omega_K \setminus \Omega_K^{(A)}$ of places at which $A$ has bad reduction (or at least a finite superset $\tilde{S}$ of $S$), and for each $v \in \Omega_K \setminus \tilde{S}$ we know how to compute the characteristic polynomial $f_v(x)$ of the Frobenius at $v$ acting on $T_\ell(A)$, where $\ell$ is any prime not divisible by $v$. These requirements are met if $A$ is the Jacobian of a genus 2 curve over $K$ given through a hyperelliptic model $y^2=f(x)$ of its affine part (here $f(x)$ is a separable polynomial of degree 5 or 6).
We construct a basic ``toolkit'' by showing that the following problems can all be solved by procedures that are guaranteed to terminate.

We pay no heed to the (astronomical) computational cost of our proposed procedures, because they will only be used to show that the problem of computing $\operatorname{End}_K(A)$ is solvable \textit{in principle}. Efficient algorithms will be discussed in section \ref{sect:ReallyDetermineEnd}.

\paragraph{Compute the characteristic polynomials of Frobenius over an extension of $K$.}
Let $K'$ be a finite extension of $K$, and let $w$ be a place of $K'$ that is above the place $v$ of $K$. Suppose that $A$ has good reduction at $v$, and let $d=[\mathbb{F}_w:\mathbb{F}_v]$ be the relative inertia degree. Then the Frobenius of $\mathbb{F}_w$ acting on $A_w$ has characteristic polynomial equal to $f_v^{[d]}(x)$, simply because the Frobenius of $\mathbb{F}_w$ is the $d$-th power of the Frobenius of $\mathbb{F}_v$.

\paragraph{Compute a number field $F \supseteq K$ over which all the endomorphisms of $A$ are defined.}\label{sect:AllEndom}
It is shown in \cite[Theorem 2.4]{MR1154704} that all the endomorphisms of $A$ are defined over $K(A[3])$. We observe that $K(A[3])$ is an extension of $K$ ramified at most at the places dividing 3 and at the places of $S$, of degree bounded by $B=|\operatorname{GL}_{2\dim A}(\mathbb{F}_3)|$. The following is thus a possible procedure to compute a suitable number field $F$: determine all extensions of $K$ of degree at most $B$ and unramified outside $S \cup \{v \in \Omega_K : p_v=3\}$ (this can be done effectively, since one of the standard proofs of Hermite's theorem is effective), and take $F$ to be the compositum of all these (finitely many) number fields. Then $K(A[3])$ is contained in $F$, hence all the endomorphisms of $A$ are defined over $F$.

Alternatively, one can also explicitly determine $K(A[3])$ by writing down polynomials whose roots are the coordinates of the $3$-torsion points of $A$ (this is possible, at least in principle, if $A$ is a genus 2 Jacobian: see the discussion in section \ref{sect:DeterminingItInPrinciple}).

\subsection{Auxiliary results about characteristic polynomials of Frobenius}

We collect in this section a few results about the sort of properties of $A$ that one might detect through the study of the associated characteristic polynomials of Frobenius.

We start with absolute irreducibility. The results of \cite{MR2496739}, specialized to the case $\dim A=2$, give:
\begin{theorem}\label{thm:Achter}
Let $A$ be an absolutely simple abelian surface such that $\operatorname{End}_{\overline{K}}(A)$ is an order in a field. The set of places $v$ of $K$ for which $A_v$ is well-defined and absolutely irreducible has density one.
\end{theorem}

We also remark that combining \cite{MR2914900} with \cite{Surfaces} (for the case of real multiplication and of trivial endomorphisms) and \cite{CM} (for the case of CM) one can give an explicit upper bound on the smallest place, as measured by its norm, for which $A_v$ is absolutely simple. In particular:
\begin{proposition}
There exists an effectively computable bound $B=B(A/K)$ such that, if $A_v$ is nonsimple for all places $v$ with $q_v \leq B$, then either $A$ itself is nonsimple or $\operatorname{End}_{\overline{K}}(A)$ is an order in a quaternion algebra.
\end{proposition}

The proof is not very different from that of proposition \ref{prop_RMBound} below, so we only sketch the argument. It suffices to show that if $\EndKBar$ is a field then there is an explicitly computable bound $B$ such that $A_v$ is irreducible for some place $v$ with $q_v \leq B$. Assuming therefore that $\EndKBar$ is a field, for $\ell$ larger than an explicit bound $\ell_0$ the image of the Galois representation is as large as it can be, given the structure of the endomorphism algebra (this is the input from \cite{Surfaces} and \cite{CM}). Suppose for simplicity that $\operatorname{End}_{\overline{K}}(A)=\mathbb{Z}$: then for any $\ell>\ell_0$ the group $G_\ell=\operatorname{GSp}_4(\mathbb{F}_\ell)$ contains some element $M$ whose characteristic polynomial is irreducible modulo $\ell$. By Chebotarev's theorem, there is a place $v \in \Omega_K^{(A)}$ such that $\rho_\ell(\operatorname{Fr}_v)=M$, and the norm of $v$ can be bounded effectively. Since the characteristic polynomial of the Frobenius at $v$ is irreducible, the reduction $A_v$ is also irreducible, and we have an effective bound on $q_v$.
For the case of real and complex multiplication (and to get much better estimates, also in the case $\operatorname{End}_{\overline{K}}(A)=\mathbb{Z}$), one obtains from \cite{Surfaces} and \cite{CM} an explicit set $\mathbb{L}$ to use as input in \cite[Lemma 4.3]{MR2914900}. Finally, notice that in order to apply the results of \cite{Surfaces} and \cite{CM} one also need to effectively bound the discriminant of the field $\EndKBar$: this can be done for example by appealing to \cite[Proposition 2.12]{GaelUniformes}.

\medskip

Next we consider the property of having a noncommutative endomorphism algebra:

\begin{lemma}\label{lemma_QMIsDetectedBySquaredPolynomials}
Let $A$ be an absolutely simple abelian surface.
The endomorphism algebra $\operatorname{End}_{\overline{K}}^0(A)$ is noncommutative (i.e. it is a division quaternion algebra) if and only if for every $v \in \Omega_K^{(A)}$ the polynomial $f_v^{[12]}(x)$ is the square of a polynomial with integral coefficients.
\end{lemma}
\begin{proof}
Suppose first that $\operatorname{End}_{\overline{K}}^0(A)$ is a division quaternion algebra.
Let $\operatorname{Fr}_v$ be a Frobenius element of $\abGal{K}$ corresponding to the place $v$. Let $\ell$ be a prime not divisible by $v$ such that $Q \otimes \mathbb{Q}_\ell$ is split (all but finitely many primes satisfy this condition). The polynomial $f_v(x)$ can be computed as the characteristic polynomial of $\rho_{\ell^\infty}(\operatorname{Fr}_v) \in G_{\ell^\infty} \subseteq \mathcal{G}_\ell$. By §\ref{sect_MT} we know that the group of connected components of $\mathcal{G}_\ell$ has exponent dividing 12, so $\rho_{\ell^\infty}(\operatorname{Fr}_v)^{12}$ belongs to 
\[
\mathcal{G}_\ell^0 = (Q \otimes \mathbb{Q}_\ell)^\times = \left\{ (x,y) \in \operatorname{End}(\mathbb{Q}_\ell^2) \oplus \operatorname{End}(\mathbb{Q}_\ell^2) \subseteq \operatorname{End}(\mathbb{Q}_\ell^4) \bigm\vert x=y \right\}.
\]
This description makes it clear that the characteristic polynomial of any element of $\mathcal{G}_\ell^0$ is the square of a polynomial with rational coefficients, and on the other hand $f_v^{[12]}(x)$ is precisely the characteristic polynomial of $\rho_{\ell^\infty}(\operatorname{Fr}_v)^{12}$ (see remark \ref{rmk_CharPoly}), so it is the square of a polynomial with integral coefficients.

For the converse implication, suppose by contradiction that $\operatorname{End}_{\overline{K}}(A)$ is not a quaternion algebra over $\mathbb{Q}$. Then the Mumford-Tate group of $A$ has rank 3 (§\ref{sect_MT}), and by a result of Serre (also independently proven by Zarhin; see \cite[Corollary 3.8]{MR1156568} for a published account of the proof) the set of places $v$ such that the eigenvalues of $\rho_{\ell^\infty}(\operatorname{Fr}_v)$ generate a rank-3 abelian subgroup of $\overline{\mathbb{Q}}^\times$ has density 1. For such a place $v$ it follows in particular that 3 of the roots of $f_v(x)$ are multiplicatively independent, hence $f_v^{[12]}(x)$ has at least 3 distinct roots, contradiction.
\end{proof}

\section{$\operatorname{End}_K(A)$ can be determined in principle}\label{sect:DeterminingItInPrinciple}
We now set out to prove that the rings $\operatorname{End}_K(A)$ and $\operatorname{End}_{\overline{K}}(A)$ can both be determined by a finite procedure. However, since this procedure is hopelessly slow, we shall then also describe more practical algorithms to determine $\operatorname{End}_{\overline{K}}(A)$ at least when $A$ is the Jacobian of a genus 2 curve over the rational numbers. Notice however that the ability to compute $\operatorname{End}_K(A)$ for any number field $A$ gives, at least in principle, a method to also determine the Galois structure of the $\operatorname{Gal}(\overline{K}/K)$-module $\operatorname{End}_{\overline{K}}(A)$, something which the knowledge of $\operatorname{End}_{\overline{K}}(A)$ alone does not give.

We start by formalizing the fact that the computation of $\operatorname{End}_K(A)$ is a strictly harder problem that the computation of $\operatorname{End}_{\overline{K}}(A)$:

\begin{lemma}
Let $g \geq 1$ Suppose that there is an algorithm to compute $\operatorname{End}_K(A)$ for any abelian variety $A$ of dimension $g$ defined over a number field $K$. Then there is an algorithm to compute $\operatorname{End}_{\overline{K}}(A)$ for any abelian variety $A$ of dimension $g$ defined over a number field $K$.
\end{lemma}
\begin{proof}
By §\ref{sect:AllEndom}, given $A/K$ there is a procedure that determines a number field $K'$, containing $K$, such that $\operatorname{End}_{K'}(A)=\operatorname{End}_{\overline{K}}(A)$. Apply the given algorithm to $A/K'$.
\end{proof}

Thus we only need to show that $\operatorname{End}_K(A)$ can be computed for any abelian surface $A$ over any number field $K$. One way to proceed is to exploit the fact that we know Tate's conjecture to hold for abelian varieties, namely, the equality
\[
\operatorname{End}_K(A) \otimes \mathbb{Z}_\ell = \operatorname{End}_{G_{\ell^\infty}} T_\ell(A)
\]
holds for all primes $\ell$.
This does not imply in general that the natural inclusion 
\[
\operatorname{End}_K(A) \otimes \mathbb{Z}/\ell^n\mathbb{Z} \subseteq \operatorname{End}_{G_{\ell^n}} A[\ell^n]
\]
is an equality when $n$ is sufficiently large; however, we have the following lemma (which applies to abelian varieties of arbitrary dimension):
\begin{lemma}\label{lemma:FiniteLevelTate}
For fixed $\ell$ and for $n$ large enough, $\operatorname{End}_K(A) \otimes \mathbb{F}_\ell$ is the reduction modulo $\ell$ of $\operatorname{End}_{G_{\ell^n}} A[\ell^n]$.
\end{lemma}
\begin{proof}
It is clear that $\operatorname{End}_K(A) \otimes \mathbb{F}_\ell$ is contained in the reduction modulo $\ell$ of $\operatorname{End}_{G_{\ell^n}} A[\ell^n]$, so we only need to show the other inclusion. Since $\operatorname{End}_{G_{\ell}} A[\ell]$ is a finite set, it suffices to show that for every $M_1 \in \operatorname{End}_{G_{\ell}} A[\ell]$ we have either $M_1 \in \operatorname{End}_K(A) \otimes \mathbb{F}_\ell$ or there exists an $n$ for which $M_1$ is not in the image of the reduction map $\operatorname{End}_{G_{\ell^n}} A[\ell^n] \to \operatorname{End}_{G_{\ell}} A[\ell]$. Yet otherwise said, it suffices to show that if a given $M_1$ is in the image of all the reduction maps $\operatorname{End}_{G_{\ell^n}} A[\ell^n] \to \operatorname{End}_{G_{\ell}} A[\ell]$, then $M_1 \in \operatorname{End}_K(A) \otimes \mathbb{F}_\ell$.
Take $M_1 \in \operatorname{End}_{G_{\ell}} A[\ell]$ and suppose that for every $n \geq 1$ we can find $M_n \in \operatorname{End}_{G_{\ell^n}} A[\ell^n]$ such that $M_n \equiv M_1 \pmod{\ell}$. For each $n\geq 1$ fix $\tilde{M}_n \in \operatorname{Aut}(T_\ell A)$ that is congruent to $M_n$ modulo $\ell^n$. By the compactness of $\operatorname{Aut}(T_\ell A)$, there is a subsequence $(\tilde{M}_{n_k})_{k \geq 1}$ of $\tilde{M}_n$ that converges to a certain $\overline{M}$ in $\operatorname{Aut} T_\ell(A)$. We can suppose without loss of generality that $\overline{M} \equiv M_{n_k} \pmod{\ell^{k}}$ for every $k$; notice that one also has $n_k \geq k$. For every $k \geq 1$ and for every $g \in G_{\ell^\infty}$ one we then obtain
\[
\begin{aligned}
g \overline{M} -\overline{M} g & \equiv g \tilde{M}_{n_k} -\tilde{M}_{n_k} g \\
& \equiv g M_{n_k} - M_{n_k} g \equiv 0 \pmod{\ell^{k}},
\end{aligned}
\]
because by assumption $M_{n_k}$ commutes with the reduction of $G_{\ell^\infty}$ modulo $\ell^{n_k}$ (and $n_k \geq k$). Since this holds for every $k$ and $n_k$ is unbounded, we have that $g$ and $\overline{M}$ commute, and since this holds for every $g$, we have that $\overline{M}$ and $G_{\ell^\infty}$ commute. By Tate's conjecture, this implies that $\overline{M}$ is in $\operatorname{End}_K(A) \otimes \mathbb{Z}_\ell$, hence $\overline{M} \bmod \ell = M_1$ is in $\operatorname{End}_K(A) \otimes \mathbb{F}_\ell$.
\end{proof}

This suggests the following (hopelessly slow) night-and-day algorithm, which we only sketch since its practical usefulness is limited. Start by embedding $A$ in some projective space $\mathbb{P}$ of dimension $N$ (this has been worked out in practice if $A$ is the Jacobian of a genus 2 curve, see \cite{MR1041476} and \cite{MR1406090}).
Then proceed as follows:
\begin{itemize}
\item by day, compute the action of Galois on $A[\ell^n]$. Since the addition law is part of the data of the abelian variety (or can be determined explicitly, if $A$ is presented as the Jacobian of a curve), one can write down polynomials whose roots are the coordinates of the $\ell^n$-division points of $A$, which in turn determines the action of Galois on them (recall that the problem of computing the Galois group of a polynomial is -- at least in principle -- effectively solved). 
From this data, compute 
\[
\operatorname{Image}(\operatorname{End}_{G_{\ell^n}} A[\ell^n] \to \operatorname{End}_{G_{\ell}} A[\ell]) \supseteq \operatorname{End}_K(A) \otimes \mathbb{F}_\ell,
\]
which gives an upper bound on $\operatorname{rank}_{\mathbb{Z}} \operatorname{End}_K(A)$.
\item by night, enumerate $(N+1)$-tuples $(P_0,\ldots,P_N)$ of homogeneous polynomials of degree at most $n$, with coefficients in $K$ and of height at most $n$, and for each, check whether the corresponding map $(P_0,\ldots,P_N) : A \darrow \mathbb{P}$ is defined over all of $A$, maps the origin of $A$ to itself, and has image contained in $A$ (in principle, these tests can all be done via Gr\"obner bases).
\end{itemize}

Lemma \ref{lemma:FiniteLevelTate} ensures that, after a finite number of steps, the endomorphisms found by night generate a $\mathbb{Z}$-module of rank equal to the upper bound found by day, so that we have at least determined the rank of $\operatorname{End}_K(A)$, and we have computed generators for an order $R$ of $\operatorname{End}_K(A)$. Now $E:=R \otimes \mathbb{Q}$ is a semisimple algebra, hence a product of matrix algebras over division rings. The index of $R$ in $\operatorname{End}_K(A)$ is certainly not larger than the index of $R$ in a maximal order of $E$, and the latter is bounded by an explicitly computable function (see \cite[Proposition A.5]{Surfaces}). Hence we can compute a finite list of primes $\ell_1, \ldots, \ell_k$ that contains all the prime divisors of $[\operatorname{End}_K(A) : R]$.

For each of these primes we now repeat the above procedure.
More precisely, we compute the image of $R$ in $\operatorname{Aut} A[\ell_i]$ (notice that this is not necessarily $R \otimes \mathbb{F}_\ell$ -- in fact, one has equality if and only if $\ell_i \nmid [\operatorname{End}_K(A) : R]$), and for increasing values of $n$ we compute $I_n:=\operatorname{Image}(\operatorname{End}_{G_{\ell_i^n}} A[\ell_i^n] \to \operatorname{End}_{G_{\ell}} A[\ell])$. If $\ell_i \nmid [\operatorname{End}_K(A) : R]$ then for $n$ large enough we will find that the image of $R$ in $\operatorname{Aut} A[\ell]$ is equal to $I_n$, which \textit{proves} $\ell_i \nmid [\operatorname{End}_K(A):R]$. Otherwise, the ``night" part of the computation will eventually find a further element of $\operatorname{End}_K(A)$, thus allowing us to enlarge $R$ (our current best guess for $\operatorname{End}_K(A)$). This process will eventually terminate with a new candidate $R$ and a certificate of the fact that $\ell_i \nmid [\operatorname{End}_K(A) : R]$. Carrying out this procedure for every $\ell_i$ finally leads to the determination of $\operatorname{End}_K(A)$.

\medskip

While the previous approach is obviously useless in practice, it suggests that Galois representations might help us determine the endomorphism ring $\operatorname{End}_K(A)$. An object which is much easier to compute than the full Galois representation attached to $A$ is the characteristic polynomial $f_v(x)$ of the Frobenius automorphism corresponding to a place $v$ of $K$ at which $A$ has good reduction: this boils down to counting points modulo $v$, and especially if $v$ has degree 1 this can be done very efficiently. On the other hand, we remark that $f_v(x)$ is an isogeny invariant, while $\operatorname{End}_K(A)$ is not: if $A_1, A_2$ are isogenous, then $\operatorname{End}_K(A_1)$ and $\operatorname{End}_K(A_2)$ are both orders in $\operatorname{End}_K(A_1) \otimes \mathbb{Q}$, so they are commensurable, but not necessarily equal. We shall show that computing characteristic polynomials of Frobenius is often enough to determine the isomorphism class of $\operatorname{End}_{\overline{K}}(A) \otimes \mathbb{Q}$, and in many cases, also finer information.

In practice, we shall only use the information coming from the characteristic polynomials to prove \textit{upper bounds} on $\operatorname{End}_{\overline{K}}(A) \otimes \mathbb{Q}$, but we remark that it is also possible to use them to certify that a given abelian variety admits non-trivial endomorphisms. In the case of abelian surfaces, this leads to an alternative procedure to obtain \textit{lower bounds} about $\EndKBar$, whose computational cost however is again prohibitive.
\begin{remark}
The advantage of working only with characteristic polynomials of Frobenius is that this does not require (even in theory) an explicit embedding of the abelian surface in projective space; in particular, one could also deal with abelian surfaces that are not Jacobians, provided that they have a way of computing Frobenius polynomials. The main drawbacks of this approach, on the other hand, are the prohibitive computational cost of establishing lower bounds and the fact that when $\EndKBar$ is a nonsplit quaternion algebra it seems very hard to determine its isomorphism class by just looking at Frobenius polynomials.
\end{remark}

We now give an indication of how the knowledge of sufficiently many Frobenius polynomials enables one to prove lower bounds on $\EndKBar$. We do not give details for all cases, but we limit ourselves to showing that computing a certain (astronomically large, but finite) number of characteristic polynomials of Frobenius allows us to certify that $A$ splits as the product of two elliptic curves (proposition \ref{prop:ProveSplitting}) or that $\operatorname{End}^0_{\overline{K}}(A)$ is a certain real number field (section \ref{sect:FindF}). With similar techniques, one can show that computing a sufficiently large (in practice, so large as to be unfeasible) number of characteristic polynomials of Frobenius, one can also prove (or disprove) the following statements: that $A$ admits (potential) complex multiplication under a given CM field $F$, that the geometric endomorphism ring of $A$ is an order in a field, and that the geometric endomorphism ring of $A$ is an order in a nonsplit quaternion algebra.

\begin{proposition}\label{prop:ProveSplitting}
There is an effective procedure which, given an abelian surface $A/K$, either finds two elliptic curves $E_1, E_2/K$ such that $A$ is $K$-isogenous to $E_1 \times E_2$, or proves that $A$ is $K$-simple. Similarly, there is a procedure that answers the same question over $\overline{K}$.
\end{proposition}
\begin{proof}
The places of bad reduction of an abelian variety are an isogeny invariant. Therefore, if there is an isogeny $A \sim E_1 \times E_2$ defined over $K$, the places of bad reduction of $E_1, E_2$ are a subset of the places of bad reduction of $A$. Let $S$ be the finite set of places of bad reduction of $A$.
One can then effectively list the (finitely many) elliptic curves over $K$ with good reduction away from $S$ (this is well-known; a computational approach to the problem is discussed for example in \cite{MR2367320}). For each pair $(E_1, E_2)$ of such elliptic curves, we can then use \cite[Lemma 1.2]{MR2181871} to check whether $A$ is isogenous to $E_1 \times E_2$. If we find such an isogeny we are done, and otherwise we have proved that $A$ is $K$-simple; this procedure obviously terminates. To get the same result over $\overline{K}$, we first find a number field $K'$ over which all the endomorphisms of $A$ are defined (see section \ref{sect:BasicAlgorithms}), and then apply the previous procedure to $A/K'$.
\end{proof}

\subsection{Proving that $A$ admits real multiplication}\label{sect:FindF}
Suppose that we have been able to prove that $\EndKBar$ is either $\mathbb{Q}$ or a real quadratic field $F$ (this is a case that happens often in practice, if one uses the algorithm we describe in section \ref{sect_Field}), and we now want to certify that there are indeed extra endomorphisms, so that $\EndKBar=F$. We describe a procedure to do so by only using characteristic polynomials of Frobenius (proposition \ref{prop_RMBound} below). We need the following preliminary lemma, which is akin to lemma \ref{lemma_QMIsDetectedBySquaredPolynomials}.
\begin{proposition}
Suppose $\EndKBar$ is a real quadratic field $F=\mathbb{Q}(\sqrt{m})$. Then for all places $v \in \Omega_K$ at which $A$ has good reduction the polynomial $f_v^{[2]}(x)$ is of the form
\begin{equation}\label{eq_SplitOverQm}
(x^2-a_vx+q_v)(x^2-\iota(a_v)x+q_v),
\end{equation}
where $a_v \in \mathcal{O}_F$ and $\iota(a_v)$ is the unique Galois conjugate of $a_v$ (or $\iota(a_v)=a_v$ if $a_v \in \mathbb{Z}$). In particular, $f_v^{[2]}(x)$ is reducible over $\mathbb{Q}(\sqrt{m})$.
\end{proposition}
\begin{remark}\label{rmk_test}
Notice that, given a characteristic polynomial of Frobenius $f_v(x)$, it is easy to check if $f_v^{[2]}(x)$ is of the form \eqref{eq_SplitOverQm}. Indeed, writing $a_v=w+z\sqrt{m}$ and $\iota(a_v)=w-z\sqrt{m}$, we find easily that $-2w$ is the coefficient of $x^3$ in $f_v^{[2]}(x)$, while $z$ satisfies $2 q_v + w^2 - m z^2=\frac{1}{2}(f_v^{[2]})''(0)$. Thus given $f_v(x)$ we can immediately find the only possible values of $w$ and $z$, plug them into the right hand side of \eqref{eq_SplitOverQm}, and check if the result agrees with the polynomial $f_v^{[2]}(x)$ we started with.
\end{remark}
\begin{proof}
Reasoning as in the proof of lemma \ref{lemma_QMIsDetectedBySquaredPolynomials}, we see that $f_v^{[2]}(x)$ is the characteristic polynomial of $\rho_{\ell^\infty}(\operatorname{Fr}_v)^2 \in \mathcal{G}_\ell^0(\mathbb{Q}_\ell) \subseteq \operatorname{GL}_2(F \otimes \mathbb{Q}_\ell)$, where the action of this last group on $T_\ell(A) \otimes \mathbb{Q}_\ell \cong \mathbb{Q}_\ell^4$ is given by the restriction from $F \otimes \mathbb{Q}_\ell$ to $\mathbb{Q}_\ell$  of the standard representation of $\operatorname{GL}_2(F \otimes \mathbb{Q}_\ell)$. In particular, $f_v^{[2]}(x)$ can be computed as $f(x)\sigma(f(x))$, where $f(x)$ is the characteristic polynomial of the action of $\rho_{\ell^\infty}(\operatorname{Fr}_v)^2 \in \operatorname{GL}_2(F \otimes \mathbb{Q}_\ell)$ and $\sigma$ is the unique nontrivial element of $\operatorname{Gal}(F/\mathbb{Q})$. The claim follows.
\end{proof}

\begin{lemma}\label{lemma_EveryReasonablePolynomialExists}
For every prime $\ell$ the finite group $\operatorname{GSp}_{4}(\mathbb{F}_\ell)$ contains elements whose characteristic polynomial $f(x)$ has the property that $f^{[2]}(x)$ is irreducible over $\mathbb{F}_\ell$. More generally, for any $n \geq 1$ there exists a bound $b(n)$ such that for all primes $\ell\geq b(n)$ the group $\operatorname{GSp}_{4}(\mathbb{F}_\ell)$ contains elements whose characteristic polynomial $f(x)$ has the property that $f^{[n]}(x)$ is irreducible over $\mathbb{F}_\ell$. One can take $b(n)=\sqrt{2n}$.
\end{lemma}
\begin{proof}
By \cite[Theorem A.1]{MR2401624} (see also \cite{MR0258855}), every reciprocal polynomial of degree 4 in $\mathbb{F}_\ell[x]$ is the characteristic polynomial of an element $g \in \operatorname{GSp}_4(\mathbb{F}_\ell)$, so it suffices to show that there exists a degree 4 reciprocal polynomial $f(x) \in \mathbb{F}_\ell[x]$ with the property that $f^{[2]}(x)$ is irreducible over $\mathbb{F}_\ell$. In order to do this, observe that the minimal polynomial of $\alpha \in \mathbb{F}_{\ell^4}$ is a degree 4 irreducible reciprocal polynomial over $\mathbb{F}_\ell$ if and only if the following two conditions are met:
\begin{itemize}
\item $\alpha \in \mathbb{F}_{\ell^4} \setminus \mathbb{F}_{\ell^2}$ -- this is equivalent to $f(x)$ being of degree 4;
\item the orbits of $\alpha$ and of $\alpha^{-1}$ under the action of $\operatorname{Gal}\left(\mathbb{F}_{\ell^4}/\mathbb{F}_{\ell^2} \right)$ coincide -- this is equivalent to $f(x)$ being reciprocal.
\end{itemize}
Choose any $\alpha \in \mathbb{F}_{\ell^4}^\times$ of exact multiplicative order $\ell^2+1$; there are $\varphi(\ell^2+1)>0$ such elements. Let $f(x)$ be its characteristic polynomial. By what we have just seen, $f(x)$ is an irreducible reciprocal polynomial of degree 4 (indeed $\alpha^{-1}=\alpha^{\ell^2}$, so $\alpha$ and $\alpha^{-1}$ have the same orbit under the action of Galois). Clearly the polynomial $f^{[2]}(x)$ vanishes on $\alpha^2$, hence to show that it is $\mathbb{F}_\ell$-irreducibile it suffices to show that $\alpha^2 \not \in \mathbb{F}_{\ell^2}$. Now if $\alpha^2$ belonged to $\mathbb{F}_{\ell^2}$ its multiplicative order (which is precisely $\frac{\ell^2+1}{(\ell-1,2)}$) would divide $\ell^2-1$, which implies $\ell^2+1 \mid 4$, a contradiction.

The proof of the more general statement is completely analogous: if $\alpha$ is an element of $\mathbb{F}_{\ell^4} \setminus \mathbb{F}_{\ell^2}$ whose minimal polynomial $f(x)$ is such that $f^{[n]}(x)$ is reducible over $\mathbb{F}_\ell$, then the same argument as above leads to the divisibility condition $\ell^2+1 \mid (\ell-1,n) (\ell^2-1) \Rightarrow \ell^2+1 \mid 2n$, which cannot happen for $\ell \geq \sqrt{2n}$.
\end{proof}

\begin{proposition}\label{prop_RMBound}
Let $F=\mathbb{Q}(\sqrt{m})$ be a real quadratic field. Assume that $\EndKBar$ is either $\mathbb{Q}$ or $F$.
There is a bound $B$, effectively computable in terms of $A$, $K$ and $F$, such that if $f_v^{[2]}(x)$ is of the form \eqref{eq_SplitOverQm} for all places $v$ with $q_v \leq B$ then $\EndKBar$ is isomorphic to $F$.
\end{proposition}
\begin{proof}
It suffices to show the following: assume that $\EndKBar=\mathbb{Q}$. Then there is an effectively computable bound $B$ such that there exists a place $v$ of good reduction for $A$ with the property that $f_v^{[2]}(x)$ is irreducible over $F$.

To show this, notice first that by \cite[Theorem 1.3]{Surfaces} there is an effectively computable bound $B_0$ such that $G_\ell=\operatorname{GSp}_{4}(\mathbb{F}_\ell)$ for all $\ell > B_0$. Fix a prime $\ell$ larger than $B_0$ and split in $\mathbb{Q}(\sqrt{m})$. Now choose an element $g$ of $G_\ell=\operatorname{GSp}_{4}(\mathbb{F}_\ell)$ whose characteristic polynomial $f(x)$ has the property that $f^{[2]}(x)$ is irreducible over $\mathbb{F}_\ell$; such an element exists by lemma \ref{lemma_EveryReasonablePolynomialExists}.

The effective Chebotarev theorem then yields the existence of a place $w$ of $K$ such that $\rho_\ell(\operatorname{Fr}_w)=g$ and $q_w$ is bounded by an effectively computable function of $A, K$, and $\ell$. Take $B$ to be equal to the bound in Chebotarev's theorem: then if by contradiction $f_v^{[2]}(x)$ were of the form \eqref{eq_SplitOverQm} for all $v$ with $q_v \leq B$, then this would in particular be true for $v=w$. It would follow that $f_w^{[2]}(x)$ splits over $\mathbb{Q}(\sqrt{m})$, hence modulo $\ell$ since $\ell$ is split in $\mathbb{Q}(\sqrt{m})$. But
\[
f_w^{[2]}(x) \equiv \text{characteristic polynomial of } \rho_\ell(\operatorname{Fr}_w) ^{[2]} \equiv f^{[2]}(x) \pmod{\ell}
\]
is irreducible modulo $\ell$ by construction, contradiction.
\end{proof}

\section{Determining $\operatorname{End}_{\overline{K}}(A)$ in practice}\label{sect:ReallyDetermineEnd}

\subsection{Outline of the algorithm}
Let $C$ be a nice genus-2 curve defined over a number field $K$ and let $J$ be its Jacobian, which is a principally polarized abelian surface over $K$. We describe a practical algorithm to (provably) determine the structure of $\operatorname{End}_{\overline{K}}(J)$, with an emphasis on the case $K=\mathbb{Q}$.
Our algorithm consists of two parts: on the one hand we want to find a ring into which $\operatorname{End}_{\overline{K}}(J)$ embeds (we call such a ring an ``upper bound" for $\operatorname{End}_{\overline{K}}(J)$), and on the other we want to show that the embedding is an isomorphism (which requires proving ``lower bounds" on $\operatorname{End}_{\overline{K}}(J)$, that is, showing that $J_{\overline{K}}$ admits sufficiently many endomorphisms). This second part has already been studied in the literature, see for example \cite{MR1648415}, \cite{MR3540944}, and \cite{MR3540958}; for a matter of completeness, however, we discuss at least some approaches to it in section \ref{sect:Certify}.
We give a rough description of our algorithm:
\begin{enumerate}
\item We try to prove that $J_{\overline{K}}$ is irreducible. This can fail for two different reasons: either because $J_{\overline{K}}$ is in fact reducible, or because it admits quaternionic multiplication. Numerical computations as in \cite{MR3540958} allow us to discriminate between these two cases, at least heuristically. Whenever we can prove geometric irreducibility by the method we describe, we also prove that $J$ does not admit quaternionic multiplication. Hence three cases arise: if we suspect that $J_{\overline{K}}$ is reducible, go to step 2; if we have proved that $J$ is geometrically irreducible, and does not have potential QM, go to step 3; if we suspect that $J$ has potential QM, go to step 4. Of course, if we cannot prove geometric irreducibility we can also run both step 2 \textit{and} step 4, and one of the two should terminate successfully and discriminate between these two cases.

For more details, section \ref{sect_Irreducibility}. 

\item If $J$ is geometrically reducible, then $C$ admits two independent maps to elliptic curves. We (try to) produce explicit expressions for these maps; 
this is sufficient to determine $\operatorname{End}_{\overline{K}}(J)$ (in fact, it is even sufficient to determine its structure as a Galois module). This has already been discussed briefly in \cite{MR3540958}, and more details are given in section \ref{sect_Reducibility}, where we also describe the computation of the maps from $C$ to elliptic curves, which was not covered in \cite{MR3540958}.

\item If $J$ is geometrically irreducible \textit{and} it does not admit quaternionic multiplication, then $\operatorname{End}_{\overline{K}}(J)$ is an order in a number field $E$, which can be either a CM field of degree 4, a real quadratic field, or (the general case) the field of rational numbers $\mathbb{Q}$. If $E$ is a CM field, and $C$ is defined over $\mathbb{Q}$, then $C$ belongs to a finite, explicit list of curves, so this case does not present any difficulties (see §\ref{sect:CM}). Otherwise we need to determine $E$; we describe an algorithm that usually does this very quickly, and which in any case computes a finite list of fields which is guaranteed to contain $E$, see section \ref{sect_Field}. This gives the desired upper bound on $\operatorname{End}_{\overline{K}}(J)$. We then discuss how to use this information to provably determine $\operatorname{End}_{\overline{K}}(J)$ (section \ref{sect:RM}).

\item If $J$ admits (potential) quaternionic multiplication, we (try to) prove this fact by showing that $J_{\overline{K}}$ admits many well-chosen real multiplications; we discuss the method in section \ref{sect:QM}. Finally, if we succeed in showing that $\EndKBar$ contains a quaternion algebra $Q$, we still need to prove that $A_{\overline{K}}$ is not the square of an elliptic curve with CM (this implies that the inclusion $Q \subseteq \EndKBar$ is in fact an equality), and we do so by using algorithm \ref{algo:ExcludeSquareCMCurve}.
\end{enumerate}

Every step of the algorithm depends on a certain bound $B$. Repeating the algorithm with increasing values of $B$ will eventually make it terminate with an explicit description of $\operatorname{End}_{\overline{K}}(J)$ and a proof that this description is correct. In practice, even with small values of $B$ we had no trouble determining $\operatorname{End}_{\overline{\mathbb{Q}}}(J)$ for the Jacobians $J$ of all the curves considered in \cite{MR3540958}, see section \ref{sect:LMFDB}.

\subsection{Irreducibility and geometric irreducibility}\label{sect_Irreducibility}

We start by remarking the following (obvious) fact:
\begin{proposition}\label{prop:NaiveIrredTest}
Suppose $A/K$ is nonsimple: then for all $v$ at which $A$ has good reduction the polynomial $f_v(x)$ is reducible in $\mathbb{Z}[x]$.
\end{proposition}

The idea is certainly not new (see for example \cite[Section 2]{MR1363577}), but it suggests the following na\"ive test to ascertain whether $A/K$ is simple: we fix a reasonable bound $B$, and compute the characteristic polynomials of Frobenius acting on $T_\ell(A)$ for all the places in $\Omega_K^{(A)}$ of norm at most $B$. Then we check whether each of these polynomials is reducible in $\mathbb{Z}$. If we find at least one characteristic polynomial that is irreducible, we have proved that $A$ is $K$-simple. Notice that this test is not guaranteed to succeed, as shown by the following example, but we shall discuss below a variant for geometrical irreducibility which can only fail in the case of potential quaternionic multiplication.

\begin{example}
One can produce examples of curves over $K$ whose Jacobian is irreducible over $K$ but not geometrically irreducible, and for which \textit{all} the characteristic polynomials $f_v(x)$ are reducible. Take for example $C:y^2=x^5-x$, considered as a curve over $K=\mathbb{Q}(i)$, and set $J:=\operatorname{Jac}(C)$. One can show that $\operatorname{End}_K(J)$ is a maximal order $R$ in the Hamilton quaternions $\mathbb{H}$, so that in particular $J$ is $K$-simple (but $C$ maps to the elliptic curve $y^2=x^3-x$, so $J$ is geometrically nonsimple). On the other hand, fix a prime $\ell \neq 2$. The Tate module $T_\ell(J)$ admits an action of $R \otimes \mathbb{Z}_\ell=\operatorname{Mat}_2(\mathbb{Z}_\ell)$, and one can show that it is a free $\operatorname{Mat}_2(\mathbb{Z}_\ell)$-module of rank 1. This implies that all the characteristic polynomials $f_v(x)$ are squares, hence a fortiori reducible.
\end{example}

There is also a variant of this test for proving $\overline{K}$-simplicity. Indeed, by the discussion in section \ref{sect_MT} we know that if $A$ is $\overline{K}$-reducible, then it is $F$-reducible for some Galois extension $F/K$ such that the inertia degree of $w/v$ divides 12 for all places $v$ of $K$ and all places $w$ of $F$ lying above $v$. Hence if $A$ is reducible over $\overline{K}$, the characteristic polynomial of $f_w(x)$ is reducible for all places $w \in \Omega_F^{(A)}$. By remark \ref{rmk_CharPoly}, we have $f_w(x)=f_v^{[d]}(x)$, where $d=f(w|v)$ is the inertia degree of $w$ over $v$. Furthermore we have $d \bigm\vert 12$ (lemma \ref{lemma:ExpDivides12}), and if $f_v^{[d]}(x)$ is reducible, then so is $f_v^{[12]}(x)$ (see remark \ref{rmk:ReducibilityTwistedPolys}). Hence:
\begin{proposition}\label{prop:GeometricIrreducibility}
Suppose $A/K$ is \textit{geometrically} reducible: then for all $v$ at which $A$ has good reduction the polynomial $f_v^{[12]}(x)$ is reducible in $\mathbb{Z}[x]$.
\end{proposition}
As with the proposition \ref{prop:NaiveIrredTest}, this leads to a simple test to try and decide whether $A_{\overline{K}}$ is irreducible.
Also notice that if $A/K$ admits QM (or $A_{\overline{K}}$ does), then $f_v(x)$ (respectively $f_v^{[12]}(x)$) is always the square of a polynomial with integer coefficients (lemma \ref{lemma_QMIsDetectedBySquaredPolynomials}), so if we find one polynomial $f_v(x)$ (resp.~$f_v^{[12]}(x)$) that is irreducible we have also proven that $A_K$ (resp.~$A_{\overline{K}}$) does not admit quaternionic multiplication.

We shall also need to deal with the particular case of $A_{\overline{K}}$ being isogenous to the square of an elliptic curve with CM. This is usually very easy to rule out:
\begin{proposition}\label{prop:NoSquareCMCurve}
Suppose $A_{\overline{K}}$ is isogenous to the square of an elliptic curve with complex multiplication by the imaginary quadratic field $F$. Then for all $v$ at which $A$ has good reduction the polynomial $f_v^{[12]}(x)$ is the square of a polynomial $g_v(x)$ in $\mathbb{Z}[x]$. Moreover, the splitting field of $g_v(x)$ is either $\mathbb{Q}$ or $F$.
\end{proposition}
\begin{proof}
Let $v \in \Omega_K^{(A)}$ and let $\mathbb{F}_v$ be the corresponding residue field. Let $K'$ be the minimal extension of $K$ over which all the endomorphisms of $A$ are defined. Fix a place $w$ of $K$ extending $v$ and let $d = [\mathbb{F}_w : \mathbb{F}_v]$; we know that $d \mid 12$. The polynomial $f_v^{[12]}(x)$ can also be obtained as $f_w^{[12/d]}(x)$, so it suffices to show the claim with $f_v^{[12]}(x)$ replaced by $f_w^{[12/d]}(x)$. But this is easy, because $A_{K'}$ is isogenous to $E^2$ for some $E/K'$ with complex multiplication defined over $K'$, hence the characteristic polynomial $f_w(x)$ is equal to the square of the characteristic polynomial of $\operatorname{Frob}_w$ acting on $E$. The splitting field of any such polynomial is either $\mathbb{Q}$ or the imaginary quadratic field $\operatorname{End}_{K'}(E)=F$ as claimed.
\end{proof}

If $A_{\overline{K}}$ is isogenous to the square of a CM elliptic curve, then none of our tests for proving ``upper bounds" gives us any information on $\operatorname{End}_{\overline{K}}(A)$, but this is to be expected because in this case $\dim_{\mathbb{Q}} \EndKBar = 8$ is as large as it can be.
Suppose now that we have certified that $\EndKBar$ contains a quaternion algebra $Q$: then 
the only two possibilities are $\EndKBar=Q$ and $A_{\overline{K}} \sim E^2$, where $E/\overline{K}$ is an elliptic curve with CM. If the equality $\EndKBar = Q$ holds, then we can prove this by using the following algorithm, whose correctness follows from proposition \ref{prop:NoSquareCMCurve}, and which terminates by proposition \ref{prop:NoSquareCMCurveTerminates} below.
\begin{algorithm}\label{algo:ExcludeSquareCMCurve} Let $A$ be an abelian surface over a number field $K$.
\begin{enumerate}
\item We loop over a positive integer $B$; given $B$, we compute $f_v^{[12]}(x)$ for all places $v \in \Omega_K^{(A)}$ of norm at most $B$.
\item For each such place $v$, we test whether $f_v^{[12]}(x)$ is a square in $\mathbb{Z}[x]$. If at least one of the polynomials $f_v^{[12]}(x)$ is not a square, then $A_{\overline{K}}$ is not isogenous to the square of an elliptic curve with CM. If they are all squares, write $f_v^{[12]}(x)=g_v(x)^2$ for some $g_v(x) \in \mathbb{Z}[x]$.
\item If $f_v^{[12]}(x)$ is a square for every $v$ of norm at most $B$, we compute the splitting field of $g_v(x)$ for every such $v$. If we find at least two nonisomorphic quadratic fields, then $A_{\overline{K}}$ is not isogenous to the square of an elliptic curve with CM. Otherwise, increase $B$.
\end{enumerate}
\end{algorithm}

If $\EndKBar$ is a quaternion algebra then this algorithm is guaranteed to terminate:
\begin{proposition}\label{prop:NoSquareCMCurveTerminates}
Suppose $A/K$ is a geometrically irreducible abelian surface such that $\EndKBar$ is a nonsplit quaternion algebra. Then algorithm \ref{algo:ExcludeSquareCMCurve} terminates.
\end{proposition}
\begin{proof}
Let $K'$ be the minimal extension of $K$ over which $A$ admits quaternionic multiplication.  
For $\ell$ large enough we have $\rho_{\ell^\infty}(\abGal{K'})=\operatorname{GL}_2(\mathbb{Z}_\ell)$, embedded block-diagonally in $\operatorname{GSp}_{4}(\mathbb{Z}_\ell)$ (see for example \cite{Surfaces}). Fix a prime $\ell > 11$ for which this equality holds. There exist $M_1, M_2 \in \operatorname{GL}_2(\mathbb{Z}_\ell)$ whose characteristic polynomials $q_1(x), q_2(x)$ have the following properties: $q_1^{[12]}(x)$ is irreducible modulo $\ell^2$; $q_2^{[12]}(x)$ is irreducible modulo $\ell$; $v_\ell(\operatorname{disc} q_1^{[12]}(x))=1$; $v_\ell(\operatorname{disc} q_2^{[12]}(x))=0$: for example, one might take $q_1(x)=x^2-2x+(1-\ell)$ and $q_2(x)=(x-g)(x-g^{\ell})$, where $g$ is a generator of $\mathbb{F}_{\ell^2}^\times$ (these polynomials have the required properties for any prime $\ell > 11$).
 Let $N_1=\operatorname{diag}(M_1,M_1)$, resp.~$N_2=\operatorname{diag}(M_2,M_2)$, be the block-diagonal operator in $\rho_{\ell^\infty}(\abGal{K'})$ whose blocks are given by $M_1$, resp.~$M_2$. 
 Since Frobenius elements are dense in $\abGal{K'}$, we can find places $w_1, w_2$ of $K'$ at which $A$ has good reduction and such that $\rho_{\ell^\infty}(\operatorname{Fr}_{w_i})$ is arbitrarily close (in the $\ell$-adic topology) to $N_i$ for $i=1,2$. In particular, by approximating $N_i$ well enough with $\rho_{\ell^\infty}(\operatorname{Fr}_{w_i})$, we can ensure that the following statements hold (we write $f_{w_1}^{[12]}(x)=g_{w_1}(x)^2, f_{w_2}^{[12]}(x)=g_{w_2}(x)^2$):
 \[
 \begin{array}{cc}
g_{w_1}(x) \text{ is irreducible modulo }\ell^2, & v_\ell(\operatorname{disc} g_{w_1}(x))=1 \\
g_{w_2}(x) \text{ is irreducible modulo }\ell, & v_\ell(\operatorname{disc} g_{w_2}(x))=0.
\end{array}
 \]
 The irreducibility conditions imply that the splitting fields of $g_{w_1}(x), g_{w_2}(x)$ are quadratic fields; the conditions on the valuations give that $\ell$ ramifies in the splitting field of $g_{w_1}(x)$ but not in that of $g_{w_2}(x)$. It now suffices to notice that algorithm \ref{algo:ExcludeSquareCMCurve} terminates as soon as it encounters both the places $v_1$ and $v_2$ of $K$ that lie respectively below $w_1, w_2$: to see this, write $d_i=[\mathbb{F}_{w_i}:\mathbb{F}_{v_i}]$ and notice that $g_{w_i}(x)=g_{v_i}^{[d_i]}(x)$. From this it follows that the at most quadratic splitting field of $g_{v_i}(x)$ contains the splitting field of $g_{w_i}(x)$, so the splitting fields of $g_{v_1}(x)$, $g_{v_2}(x)$ are respectively the same as the splitting fields of $g_{w_1}(x)$, $g_{w_2}(x)$, which we have shown not to be isomorphic.
\end{proof}

\medskip

With an eye to practical implementation, we also describe some additional tests one can use to prove that $\operatorname{End}_{\overline{K}}(A)=\mathbb{Z}$ when this is the case. We start with the following fact, which does not seem to have been observed before and which can be useful in practice:
\begin{proposition} Let $f(x) \in K[x]$ be an irreducible polynomial of degree 5. Let $C$ be the smooth irreducible curve corresponding to the affine model $y^2=f(x)$. The abelian variety $\operatorname{Jac}(C)$ is absolutely irreducible.
\end{proposition}
\begin{proof}
We shall make use of the following facts, which are all easy to check:
\begin{enumerate}
\item $\operatorname{GSp}_{4}(\mathbb{F}_2)$ is isomorphic to $S_6$; in particular, all elements of order 5 in $\operatorname{GSp}_{4}(\mathbb{F}_2)$ are conjugate to each other;
\item all transitive subgroups of $S_5$ have order divisible by 5;
\item the centralizer in $\operatorname{End}(\mathbb{F}_2^4)$ of any element $g \in \operatorname{GSp}_4(\mathbb{F}_2)$ of exact order 5 consists of precisely $16$ elements, which are all invertible with the only exception of the zero matrix (using fact 1, it is enough to check this for a single element of order 5, and the verification is immediate. We remark that an element of order 5 in $\operatorname{GSp}_4(\mathbb{F}_2)$ is given by $\begin{pmatrix}
1 & 0 & 0 & 1 \\0 & 0 & 1 & 1\\0 & 1 & 1 & 1 \\1 & 1 & 1 & 1
\end{pmatrix}$);
\item if $\operatorname{Jac}(C)$ is nonsimple over an extension $F$ of $K$, then $\operatorname{End}_{F}(\operatorname{Jac}(C)) \otimes_{\mathbb{Z}} \mathbb{F}_2$ contains nontrival idempotents. 
\end{enumerate}
We now prove that $J:=\operatorname{Jac}(C)$ is irreducible over $K$. Indeed, suppose that $J$ is reducible: then $\operatorname{End}_{K}(J) \otimes_{\mathbb{Z}} \mathbb{F}_2$ contains a nontrivial idempotent, which is not an invertible operator when considered as an endomorphism of $J[2] \cong \mathbb{F}_2^4$. The elements of $\operatorname{End}_{K}(J) \otimes_{\mathbb{Z}} \mathbb{F}_2$ commute with the action of Galois on $J[2]$. 

On the other hand, let $d_1,\ldots,d_5$ be the roots in $\overline{\mathbb{Q}}$ of the equation $f(x)=0$, and let $P_i=(d_i,0) \in C(\overline{\mathbb{Q}})$. It is well-known that $J[2]$ is the $\mathbb{F}_2$-vector space generated by the classes of the divisors $d_i:=(P_i)-(0)$, subject to the only relation $d_1+\ldots+d_5=0$.
It follows easily that $J[2]$ is a faithful representation of $\operatorname{Gal}(f(x))$.

Now observe that $\operatorname{Gal}(f(x))$ contains an element of order 5 (if $f(x)$ is irreducible of degree 5, its Galois group is a transitive subgroup of $S_5$), so the elements of $\operatorname{End}_{K}(J) \otimes_{\mathbb{Z}} \mathbb{F}_2$, including the nontrivial idempotent, commute with an element of order 5 in $\operatorname{GSp}_4(\mathbb{F}_2)$. But this contradicts fact (3).

On the other hand, if $F/K$ is a normal extension of $K$ such that $5 \nmid [F:K]$, then the Galois group of (the splitting field of) $f(x) \in F[x]$ still has order divisible by 5, hence the same argument applies to show that $J_F$ is simple over $F$. Finally, if $J_{\overline{K}}$ is nonsimple, then all the endomorphisms of $J$ are defined over a solvable normal extension $F$ of $K$ of degree dividing 48 (see section \ref{sect_MT}), and by what we have already seen this finishes the proof.
\end{proof}

Finally, the following two criteria (the former of which is a theorem of Zarhin) can be used to speed up the process of proving that $\operatorname{End}_{\overline{K}}(J)=\mathbb{Z}$ when this is the case.
\begin{theorem}
Let $C$ be the smooth projective curve (of genus at least 2) corresponding to the affine model $y^2=f(x)$. Suppose that either of the following two conditions is satisfied:
\begin{enumerate}
\item (Zarhin) the Galois group of $f(x)$ is $A_n$ or $S_n$, where $n=\deg(f(x))$;
\item $\operatorname{Jac}(C)$ is geometrically simple, and there exists a place $v$ of $K$ such that the $v$-valuation of $\operatorname{disc}(f(x))$ is equal to 1.
\end{enumerate}
Then $\operatorname{End}_{\overline{K}}(\operatorname{Jac}(C))=\mathbb{Z}$.
\end{theorem}
\begin{remark}
The second condition is inspired by Hall's paper \cite{MR2820155}.
\end{remark}

\begin{proof}
The first statement is the main result of \cite{MR1748293}.
As for the second, the hypothesis that $v(f(x))=1$ implies (see \cite[§9.2, Example 8]{MR1045822} or the appendix to \cite{MR2820155}) that $A$ is semistable at $v$, of toric rank 1. The toric rank of a semistable abelian variety is invariant under finite extension of the ground field, so we can assume that $\operatorname{End}_K(A)=\operatorname{End}_{\overline{K}}(A)$ and that there is a place $\mathfrak{p}$ of $K$ at which the special fiber of the Néron model $\mathcal{A}$ of $A$ has toric dimension 1. By definition, this means that $\mathcal{A}_\mathfrak{p} \times_{\mathcal{O}_{\mathfrak{p}}} \mathbb{F}_{\mathfrak{p}}$ is the extension of an abelian variety by a torus $T$ of rank 1. Let $R:=\operatorname{End}_{\overline{K}}(A)$, let $X$ be the character group of $T$, and write $X_\mathbb{Q}$ for $X \otimes \mathbb{Q}$. The $\mathbb{Q}$-vector space $X_{\mathbb{Q}}$ is an $(R \otimes \mathbb{Q})$-module, and since by assumption $A$ is absolutely simple we know that $\operatorname{End}^0(A)=R \otimes \mathbb{Q}$ is a division algebra. In particular, the action of $\operatorname{End}^0(A)$ on $X_\mathbb{Q}$ is either trivial or faithful, and the first case is impossible since $-1 \in \operatorname{End}^0(A)$ acts as $-1$ on $X_{\mathbb{Q}}$. It follows that $X_\mathbb{Q} \cong \mathbb{Q}$ is a faithful $\operatorname{End}^0(A)$-module, hence $\dim_{\mathbb{Q}} X_\mathbb{Q} =1$ is divisible by $\operatorname{dim}_{\mathbb{Q}}(\operatorname{End}^0(A))$: this forces $\operatorname{End}^0(A)=\mathbb{Q}$ and $R=\mathbb{Z}$.

\end{proof}

Thus our proposed algorithm to test for geometric irreducibility is as follows:
\begin{algorithm}\label{algo:Irreducibility}
Let $A/K$ be an abelian surface, presented as the Jacobian of a genus-2 hyperelliptic curve $C: y^2=f(x)$. Fix a bound $B$.
\begin{enumerate}
\item If $f(x)$ is irreducible of degree 5, output ``$A$ is absolutely irreducible" and terminate the algorithm.
\item Compute the Galois group of $f(x)$. If it is isomorphic to either $A_n$ or $S_n$ (where $n:=\deg f(x)$), return ``$\operatorname{End}_{\overline{K}}(\operatorname{Jac}(C))=\mathbb{Z}$. $A$ is absolutely irreducible and does not admit potential QM".
\item Compute $f_v^{[12]}(x)$ for all $v \leq B$. If at least one such polynomial is irreducible, return ``$A$ is absolutely irreducible" and terminate. If at least one such polynomial is not the square of a polynomial with rational coefficients, output ``$A$ does not admit potential QM" and terminate.
\item Return ``This test cannot prove that $A$ is absolutely irreducible, nor that it does not admit potential QM".
\end{enumerate}
\end{algorithm}

To obtain the version for $K$-irreducibility, simply replace $f_v^{[12]}(x)$ by $f_v(x)$ in the above. 
When $A$ is geometrically irreducible, we expect this test to be able to quickly prove that this is the case, unless $A$ admits quaternionic multiplication (see lemma \ref{lemma_QMIsDetectedBySquaredPolynomials}). In fact, theorem \ref{thm:Achter} guarantees that testing reducibility of characteristic polynomials of Frobenius \textit{will} eventually prove irreducibility if we check a sufficiently large number of places -- again, with the only exception of abelian surfaces with quaternionic multiplication. Furthermore, since the places for which $f_v^{[12]}(x)$ is irreducible have density 1 when $A/K$ is absolutely irreducible without potential QM, we expect this test to only need a fairly small number of places.

This expectation turns out to be correct for curves over $\mathbb{Q}$: in all the cases we tested (see section \ref{sect:NumericalResults}), for Jacobians of curves defined over the rationals setting $B=59$ was sufficient to show that $A$ was geometrically irreducible in all cases in which this was true and $A$ did not admit potential QM.

\subsection{Reducibility}\label{sect_Reducibility}
Suppose that we suspect (based on the previous test, or on numerical computations like those described in \cite{MR3540958}) that $J$ is geometrically reducible. We can try to prove this by finding two elliptic curves $E_1, E_2$ and two maps $\varphi_i : C \to E_i$ such that the pullbacks $\varphi_i^* \omega_i$ of the canonical differentials on $E_1, E_2$ are linearly independent. This is clearly enough to prove that $J$ is (geometrically) isogenous to the product $E_1 \times E_2$, and since the maps will be found explicitly, one can also determine the minimal field over which the splitting happens. 

An algorithm to carry out such computations has been sketched in \cite{MR3540958}; we repeat it briefly here for the convienience of the reader, and we add some remarks on the computation of the maps $\varphi_i$, which was not described in \cite{MR3540958}. We use MAGMA as our computational system of reference.

Starting from the equation of the curve $C: y^2=f(x)$, $f(x) \in K[x]$, we can compute to high numerical precision the period matrix $\Pi \in \operatorname{M}_{2,4}(\mathbb{C})$ of $J$. We then look for pairs of matrices $(A,M)$ with $A \in \operatorname{M}_2(\mathbb{C})$ and $M \in \operatorname{M}_4(\mathbb{Z})$ such that $A\Pi = \Pi M$ up to a very small numerical error. Here $A$ is the so-called analytic representation of an endomorphism, namely a map
$
\mathbb{C}^2 \to \mathbb{C}^2
$
that induces, passing to the quotient by the lattice $\Lambda$ generated by the columns of $\Pi$, a map
\[
J = \mathbb{C}^2/\Lambda \to \mathbb{C}^2/\Lambda = J
\]
that is an endomorphism of $J$. Using numerical linear algebra one can find such matrices $M$ with exact integral entries, and numerical approximations to $A$. Using standard algorithms for the recognition of algebraic numbers, we can then represent the entries of $A$ as exact algebraic numbers, which allows us to determine a number field $F$ which contains all the coefficients of the matrices $A$. The set of such matrices $A$ is a $\mathbb{Z}$-algebra $\mathcal{E}$, and one can efficiently find idempotents in $\mathcal{E} \otimes \mathbb{Q}$. With the exception of the identity and of the null matrix, a suitable multiple of such an idempotent corresponds to a pair $(A,M)$ which furthermore satisfies $\operatorname{rk}(M)=2$. We then suspect that $J$ is in fact isogenous to $E_1 \times E_2$ for two elliptic curves $E_1, E_2$, and that $A \in \operatorname{End}(J)$ corresponds to the projection on one of the two factors $E_i$, say $E_1$ (it is of course possible that $E_1, E_2$ are in turn isogenous, but this is not important for the purposes of the algorithm).

Choose two columns of $M$ that generate the (rank 2) image of $M$, call them $C_1$ and $C_2$, and fix a vector $v=(v_1,v_2) \in \mathbb{Z}^2$. For most choices of $v$, the linear map $\pi : \mathbb{C}^2 \to \mathbb{C}$ given by $(z_1,z_2) \mapsto v_1z_1+v_2z_2$ will send $\Pi C_1, \Pi C_2$ to two $\mathbb{Z}$-linearly independent complex numbers, which we will denote by $\omega_1,\omega_2$ (in fact, one can always take $v$ to be one of $(0,1)$, $(1,0)$, or $(1,1)$).

Denote by $\Lambda_E$ the lattice generated in $\mathbb{C}$ by the images of $\Pi C_1, \Pi C_2$ via $\pi$; the complex torus $\mathbb{C}/\Lambda_E$ should be an elliptic curve $E$ that $C$ maps to. We compute numerically the $j$-invariant of $\mathbb{C}/\Lambda_E$ (and recognize it as an algebraic number), and, using the classical elliptic functions $g_4(\Lambda_E), g_6(\Lambda_E)$, obtain an equation for $E$ over its minimal field of definition.

We now consider the (analytic) Jacobian of $E$, which comes with an associated period matrix $P_2$. We can numerically compute a complex number $\alpha_E$ and a matrix $M_E \in \operatorname{M}_2(\mathbb{Z})$ such that
\[
\alpha_E \begin{pmatrix}
\omega_1 \bigm\vert \omega_2
\end{pmatrix} = P_2 M_E.
\]
The composition $\psi:=\alpha_E \circ \pi \circ A$ induces a map $\mathbb{C}^2 / \Lambda \to \mathbb{C}/\Lambda_E$ which should correspond to some map $\varphi : C \to E$. Our purpose is to determine $\varphi$.

We first determine the degree of $\varphi$. Notice that $\psi : \mathbb{C}^2 \to \mathbb{C}$ has kernel isomorphic to $\mathbb{C}$, and furthermore it induces a map $\Lambda \rightarrow \Lambda_E$. Passing to the quotient by $\ker \psi$, we obtain a map $\mathbb{C} \to \mathbb{C}$, a lattice $\tilde{\Lambda}$ in $\mathbb{C}$, and a map (which we still denote by $\psi$) from $\tilde{\Lambda}$ to $\Lambda_E$. The lattice $\tilde{\Lambda}$ has rank 2, that is, it is a full lattice in $\mathbb{C}$. Indeed, denoting by $T_i$ the $i$-th column of the matrix $T$, one has
\[
\begin{aligned}
\ker \psi \cap \Lambda & = \ker(\alpha \pi A) \cap \Lambda \\
& = \left\{ \sum_{i=1}^4 \alpha_i\Pi_i : \alpha_i \in \mathbb{Z}, \; \sum_{i=1}^4 \alpha_i \pi A \Pi_i =0 \right\} \\
& = \left\{ \sum_{i=1}^4 \alpha_i\Pi_i : \alpha_i \in \mathbb{Z}, \; \sum_{i=1}^4 \alpha_i \pi (\Pi M)_i =0 \right\}:
\end{aligned}
\]
since the column space of $M$ has rank 2, this $\mathbb{Z}$-module has also rank 2, hence the kernel of $\Lambda \to \tilde{\Lambda}$ has rank 2 and $\tilde{\Lambda}$ is free of rank 2. We now have an explicit map $\psi: \mathbb{C}/\tilde{\Lambda} \to \mathbb{C}/\Lambda_E$, and we can compute $\psi^{-1}(\Lambda_E)$. The order of the finite quotient group $\psi^{-1}(\Lambda_E) \bigm / \tilde{\Lambda}$ is the order $d$ of $\varphi$.

Now we turn to the computation of equations for $\varphi$. Recall that $C$ is given by the equation $y^2=f(x)$, where $f(x) \in K[x]$.
Given any $n \in \mathbb{Z}$, we now have a point $Q_n:=(n,\sqrt{p(n)})$ on $C$, which we can map to the AnalyticJacobian of $C$, then to $\mathbb{C}/P_2$ using $\alpha_E \circ \pi \circ A$, and then back to $E$. The resulting point should be a point of $E$ defined over the same field as $Q_n$, namely $K(\sqrt{p(n)})$. We can therefore employ an LLL-based algorithm to detect $K$-linear relations between $1$, $\sqrt{p(n)}$, and the coordinates of $\varphi(Q_n)$, and therefore guess the exact value of $\varphi(Q_n)$. Write $\varphi$ as a pair $(w(x,y), z(x,y))$ of rational functions; since $y^2$ is a function of $x$, we can write
\[
w(x,y)=\frac{\sum_{i=0}^e a_i x^i + y \sum_{i=0}^{e} b_i x^i}{\sum_{i=0}^e c_i x^i}
\]
for some $e$, and similarly for $z(x,y)$. We can assume without loss of generality that the chosen model for the elliptic curve $E$ is purely hyperelliptic; under this assumption, $w(x,y)$ can be taken to be a function of $x$ alone. Indeed, the hyperelliptic involution on $C$ induces multiplication by $-1$ on $\operatorname{Jac}(C)$, and if $\psi : \operatorname{Jac}(C) \to E$ is a map of abelian variety we have $[-1]_E \circ \psi = \psi \circ [-1]_C$. If $E$ is expressed by a purely hyperelliptic Weierstrass model $z^2=g(w)$, the map $[-1]_E$ is given by $(w,z) \mapsto (w,-z)$. It follows that $w(x,-y)=w(x,y)$, so that $w$ (which is of degree at most 1 in $y$) is actually independent of $y$; we write $w(x,y)$ simply as $w(x)$. Likewise, $z(x,y)/y$ is a function of $x$ alone. In this situation the degree of $w(x)$ coincides with that of $\varphi$, which the numerical computations suggests to be $d$, so we can write
$
\displaystyle w(x)=\frac{\sum_{i=0}^d a_i x^i}{\sum_{i=0}^d c_i x^i}.
$
Using the strategy above to find the exact value of $\varphi(Q_n)$ for many different integers $n$, we now have a large number of linear equations for the unknowns $a_i, c_i$, which can be easily solved. Once an expression for $w$ is known it is a trivial matter to deduce a corresponding expression for $z$; moreover, we can then simply \textit{check} whether $(w,z)$ really gives a morphism $C \to E$, in which case we have proved that $E$ is a quotient of $C$ of degree $d$, hence that $\operatorname{Jac}(C)$ splits -- over an explicit number field -- as the product of two elliptic curves.

Using the same procedure on the complementary idempotent in $\mathcal{E} \otimes \mathbb{Q}$, we obtain a pair of explicit maps $\varphi_i : C \to E_i$. Let $\tau_i$ be the canonical differential of $E_i$: we can now compute $\nu_i:=\varphi_i^* \tau_i$ for $i=1,2$, and if we find that $\nu_1, \nu_2$ are linearly independent we have proved that $(\varphi_1,\varphi_2)$ induces an isogeny $\operatorname{Jac}(C) \sim E_1 \times E_2$. Since the maps are explicit, it is now an easy matter to compute $\operatorname{End}(\operatorname{Jac}(C))$, even as a Galois module.

To prove that we have really determined $\operatorname{End}_{\overline{K}}(\operatorname{Jac}(C))$, it actually remains to show that we have found the isogeny of minimal degree: this can be shown by verifying that $\operatorname{Jac}(C)$ belongs to the Humbert surface $H_{d^2}$, see section \ref{sect:RM} below.
We conclude this paragraph by remarking that if we are not interested in the curves $E_1, E_2$, but only in proving that $J_{\overline{K}}$ is reducible, then we can also just test whether the corresponding point in moduli space lies on a Humbert surface of square discriminant.

\begin{example}
Even for bielliptic curves, the equations of the covering can be fairly complicated. Consider for example the curve $C:y^2=f(x)=x^6+2x^5+7x^4+8x^3+11x^2+6x+5$ (\cite[\href{http://www.lmfdb.org/Genus2Curve/Q/1088/a/1088/1}{Curve 1088.a.1088.1}]{lmfdb}). According to the LMFDB, $J=\operatorname{Jac}(C)$ splits as the product of two non-isogenous elliptic curves over $\mathbb{Q}(\sqrt{2})$. To prove this fact, we look for an explicit map $\varphi:C \to E$, where $E$ is an elliptic curve over $\mathbb{Q}(\sqrt{2})$. Using MAGMA's AnalyticEndomorphisms, we see that the endomorphism ring of $J$ contains $M=\begin{pmatrix}
-1 & 0 & 0 & 1 \\
0 & -1 & -1 & 0 \\
0 & -1 & -1 & 0 \\
1 & 0 & 0 & -1
\end{pmatrix}$, of rank 2. We recognize the corresponding $A$-matrix to be $\begin{pmatrix}
-1-\frac{\sqrt{2}}{2} & -\frac{\sqrt{2}}{2} \\
- \frac{\sqrt{2}}{2} & \frac{\sqrt{2}-2}{2}
\end{pmatrix}$,
which is indeed defined over $\mathbb{Q}(\sqrt{2})$. Following the procedure just outlined, we determine the corresponding elliptic curve quotient, which is $z^2=w^3+\frac{1}{12}(7-6\sqrt{2})+\frac{1}{108}(29\sqrt{2}-36)$, and we compute the image of $(n,\sqrt{f(n)})$ for $n=-2,\ldots,2$ through our putative map $\varphi$. Using an LLL-type algorithm to detect linear relations between $1, \sqrt{2}$ and the $w$-coordinate of $\varphi$, we find
\begin{center}
\begin{tabular}{|c|c|}
\hline $n$ & $w(\varphi(n,\sqrt{f(n)}))$ \\ \hline
$-2$ & $\frac{1}{6}(345-242\sqrt{2})$ \\
$-1$ & $\frac{1}{6}(15-8\sqrt{2})$ \\
$0$ & $\frac{1}{6}(21+10\sqrt{2})$ \\
$1$ & $\frac{1}{6}(21+10\sqrt{2})$\\
$2$ &  $\frac{1}{294}(465+22\sqrt{2})$\\ \hline
\end{tabular}
\end{center}
This quickly leads us to the following expression for $w(\varphi(x,y))$:
\[
w(\varphi(x,y))=\frac{(15-8\sqrt{2})((153+4\sqrt{2})+4(14+\sqrt{2})x+97x^2)}{6 \cdot 97 (1-\sqrt{2}+x)^2},
\]
whose coefficients have quite large height compared with those of $f(x)$. The formula for $z(\varphi(x,y))$ is simpler: one has 
$
z(\varphi(x,y))=\frac{-9+7\sqrt{2}}{2(1-\sqrt{2}+x)^3}y.
$ Finally, by Poincaré's complete reducibility theorem we know that $J \sim E \times E_2$ for some other elliptic curve $E_2/\mathbb{Q}(\sqrt{2})$, and it is easy to see that $E_2$ must be the Galois conjugate ${}^\sigma E$ of $E$ (indeed, given a map $\varphi:C \to E$, the Galois conjugate of $\varphi$ is a map $C \to {}^\sigma E$).
\end{example}

\subsection{Computing $\operatorname{End}^0_{\overline{K}}(A)$ when it is a field}\label{sect_Field}
Let $A/K$ be an abelian surface such that $\operatorname{End}^0_{\overline{K}}(A)$ is a number field $E$. We give an algorithm that attempts to determine $E$ by computing its discriminant.

 By the discussion in section \ref{sect_MT} (see in particular lemma \ref{lemma:ExpDivides12}), we know that there exists an extension $K'$ of $K$ over which all the endomorphisms of $A$ are defined and which satisfies $[K':K] \mid 4$. Let $w$ be a place of $K'$ lying over a place $v$ of $K$. The degree of the extension $\mathbb{F}_w/\mathbb{F}_v$ divides 4, so the characteristic polynomial $f_w(x)$ is of the form $f_v^{[d]}(x)$ for some $d \in \{1,2,4\}$. By remark \ref{rmk:ReducibilityTwistedPolys}, if $f_v^{[4]}(x)$ is irreducible then $f_v^{[d]}(x)$ is also irreducible for $d=1,2$, so that in particular we deduce:
\begin{lemma}
Let $A/K$ be an abelian surface such that $\operatorname{End}^0_{\overline{K}}(A)$ is a number field $E$. Let $K'$ be the minimal extension over which all the endomorphisms of $A$ are defined. If $v \in \Omega_K^{(A)}$ is such that $f_v^{[4]}(x)$ is irreducible, then for all places $w$ of $K'$ lying over $v$ the abelian variety $A_w$ is irreducible.
\end{lemma}

This leads to introducing the set
\begin{equation}\label{eq:OmegaPrime}
\EPlaces:=\left\{v \in \Omega_K^{(A)} : \begin{array}{c}
f_v^{[4]}(x) \text{ is irreducible} \\
f_v(x)=x^4+ax^3+bx^2+q_vax+q_v^2, p_v \nmid b
\end{array}\right\}.
\end{equation}
Notice that as long as one can compute $f_v(x)$ it is a simple matter to test whether a place of $K$ belongs to $\EPlaces$.
By theorem \ref{thm:Achter} we know that the set of places for which $A_v$ is absolutely irreducible has density 1, so by proposition \ref{prop:DensityOrdinaryReduction} we deduce that $\Omega_K'$ has positive density (in fact, density equal to either $1$, $1/2$ or $1/4$).

\begin{theorem}\label{thm:DetermineF}
Let $A/K$ be an abelian surface such that $\operatorname{End}^0_{\overline{K}}(A)$ is a field $E$. 
For every $v \in \EPlaces$ let $F(v)$ be the CM field generated by $\pi_v$, where $\pi_v$ is the Frobenius automorphism of $A_v$ (concretely, $F(v)$ is the field generated over $\mathbb{Q}$ by one root of the polynomial $f_v(x)$). For every subset $T$ of $\EPlaces$ we have
\begin{equation}\label{eq_discDivisibility}
\operatorname{disc}(E) \bigm\vert \operatorname{gcd}_{v \in T} \operatorname{disc}(F(v)).
\end{equation}
More precisely, if $\EndKBar$ is either $\mathbb{Q}$ or a quartic CM field we have
\begin{equation}\label{eq_DiscOtherFields}
\operatorname{disc}(E)=\operatorname{gcd}_{v \in \EPlaces} \operatorname{disc}(F(v)),
\end{equation}
while when $\EndKBar$ is a real quadratic field $E$ we have
\begin{equation}\label{eq_DiscQuadField}
\operatorname{disc}(E)^2=\operatorname{gcd}_{v \in \EPlaces} \operatorname{disc}(F(v)).
\end{equation}
\end{theorem}

For the proof we need a lemma:
\begin{lemma}\label{lemma_ANT1}
Let $E/F/\mathbb{Q}$ be a tower of extensions of number fields with $[F:\mathbb{Q}]=[E:F]=2$. Let $p$ be a prime and let $v$ be a place of $F$ lying over $p$. If $\mathcal{O}_E \otimes \mathbb{Z}_p \cong (\mathcal{O}_F \otimes \mathbb{Z}_p)^2$ then the place $v$ is unramified in the relative extension $E/F$.
\end{lemma}
\begin{proof}
Suppose by contradiction that $v$ ramifies in $E$. Then $\operatorname{Spec}\left(\mathcal{O}_E \otimes \mathbb{Z}_p \right)$ has as many connected components as $\operatorname{Spec}\left(\mathcal{O}_F \otimes \mathbb{Z}_p \right)$, but this contradicts the hypothesis.
\end{proof}

We can now prove theorem \ref{thm:DetermineF}.

\begin{proof}
We start by proving \eqref{eq_discDivisibility}. 
Let $K'=\Kconn$ be the minimal extension of $K$ over which all the endomorphisms of $A$ are defined, and let $w$ be a place of $K'$ above $v$. 
Using corollary \ref{cor_ordImpliesEndDefBase} we find
\[
E=\operatorname{End}^0_{K'}(A) \hookrightarrow \operatorname{End}^0_{\mathbb{F}_w}(A_w)=\operatorname{End}^0_{\mathbb{F}_v}(A_v)=F(v),
\]
so in particular $\operatorname{disc}(E)$ divides $\operatorname{disc}(F(v))$. Notice that if $E$ is a quartic CM field then the previous injection is in fact an equality, because $F(v)$ is also a quartic number field by definition.
We now prove \eqref{eq_DiscOtherFields} and \eqref{eq_DiscQuadField} as a consequence of theorem \ref{thm_Zarhin}.
We distinguish three cases:
\begin{itemize}
\item $E$ is a quartic CM field: we have already seen that the equality $E=F(v)$ holds for all $v \in \EPlaces$, so there is nothing to prove.
\item $E=\mathbb{Q}$. It suffices to prove that for every prime $p$ there exists a place $v \in \EPlaces$ such that $p \nmid \operatorname{disc} F(v)=\operatorname{disc} \operatorname{End}^0_{\mathbb{F}_v}(A_v)$.
We know from §\ref{sect_MT} that $\mathcal{G}_\ell=\operatorname{GSp}_{4,\mathbb{Q}_\ell}$ and $K^{\text{conn}}=K$; choose a basis of $V_\ell(A)$ in such a way that the symplectic form preserved by $\mathcal{G}_\ell$ (up to similitudes) is the one associated with the matrix $
\begin{pmatrix}
0 & 1 & 0 & 0 \\
-1 & 0 & 0 &0 \\
0 & 0 & 0 & 1 \\
0 & 0 & -1 & 0
\end{pmatrix}
$. 
 Suppose for the moment that $p \neq 2$. Fix an auxiliary integer $q \neq p$ that is not congruent to $0$ or $-1$ modulo $p$ (for example $q=p+1$), and set
\[
f_p=\left(
\begin{array}{cccc}
 0 & q & 0 & 0 \\
 1 & 0 & 0 & 0 \\
 0 & 0 & 1 & q+1 \\
 0 & 0 & 1 & 1 \\
\end{array}
\right) \in \mathcal{G}_p(\mathbb{Q}_p).
\]
The roots $\pm \sqrt{q}, 1 \pm \sqrt{1+q}$ of the characteristic polynomial of $f_p$ are all distinct. Furthermore, an immediate computation shows that the centralizer of $f_p$ in $\operatorname{End} T_p(A)$ is the set of matrices of the form
\[
\left(
\begin{array}{cccc}
 a & bq & 0 & 0 \\
 b & a & 0 & 0 \\
 0 & 0 & c & d(q+1) \\
 0 & 0 & d & c \\
\end{array}
\right),
\]
which in turn is isomorphic to $\mathbb{Z}_p[\sqrt{q}] \oplus \mathbb{Z}_p[\sqrt{1+q}]$ as a ring. By theorem \ref{thm_Zarhin}, we know that there is a positive-density set of places of $K$ such that $\operatorname{End}_{\mathbb{F}_v}(A_v) \otimes \mathbb{Z}_p \cong \mathbb{Z}_p[\sqrt{q}] \oplus \mathbb{Z}_p[\sqrt{1+q}]$. By theorem \ref{thm:Achter} and proposition \ref{prop:DensityOrdinaryReduction}, there is also a positive-density subset of places $v \in \EPlaces$ such that $\operatorname{End}_{\mathbb{F}_v}(A_v) \otimes \mathbb{Z}_p \cong \mathbb{Z}_p[\sqrt{q}] \oplus \mathbb{Z}_p[\sqrt{1+q}]$. Consider such a $v$. The ring $\operatorname{End}_{\mathbb{F}_v}(A_v)$ is contained in $\mathcal{O}_{F(v)}$, and since 
$\operatorname{End}_{\mathbb{F}_v}(A_v) \otimes \mathbb{Z}_p \cong \mathbb{Z}_p[\sqrt{q}] \oplus \mathbb{Z}_p[\sqrt{1+q}]$ is integrally closed we must have $\operatorname{End}_{\mathbb{F}_v}(A_v) \otimes \mathbb{Z}_p \cong \mathcal{O}_{F(v)} \otimes \mathbb{Z}_p$. Since $p$ is unramified in $\mathbb{Z}_p[\sqrt{q}] \oplus \mathbb{Z}_p[\sqrt{1+q}]$, it follows that it is also unramified in $\mathcal{O}_{F(v)}$. This shows that, as claimed, for $p \neq 2$ there is a positive-density subset of $\EPlaces$ such that $p \nmid \operatorname{disc} F(v)$.

The argument for $p=2$ is completely analogous, the only difficulty being that $2$ ramifies in all rings of the form $\mathbb{Z}_2[\sqrt{q}]$. This is solved by taking for example $f_2=\left(
\begin{array}{cccc}
 1 & -2 & 0 & 0 \\
 2 & 3 & 0 & 0 \\
 0 & 0 & 2 & -1 \\
 0 & 0 & 1 & 3 \\
\end{array}
\right) \in \mathcal{G}_2(\mathbb{Q}_2)$, whose centralizer is isomorphic to $\mathcal{O}^2$, where $\mathcal{O}=\mathbb{Z}_2\left[ \frac{1+\sqrt{-3}}{2} \right]$ is the ring of integers of the unique unramified quadratic extension of $\mathbb{Q}_2$.

\item $E$ is a real quadratic field. By §\ref{sect_MT} we know that 
\[
\mathcal{G}_\ell^0=\mathbb{G}_m \cdot \operatorname{Res}_{E \otimes \mathbb{Q}_\ell/\mathbb{Q}_\ell}(\operatorname{SL}_{2,E \otimes \mathbb{Q}_\ell}),
\]
or more explicitly 
\[
\mathcal{G}^0_\ell(\mathbb{Q}_\ell) \cong \{ x \in \operatorname{GL}_2(E \otimes \mathbb{Q}_\ell) \bigm\vert \det(x) \in \mathbb{Q}_\ell^\times \}.
\]

Let $v$ be any place in $\EPlaces$. The quartic number field
\[
F(v) =\operatorname{End}^0_{\mathbb{F}_v}(A_v) = \operatorname{End}^0_{\overline{\mathbb{F}_v}}(A_v)
\]
contains $\EndKBar=E$, so its discriminant is
$
N_{E/\mathbb{Q}}(\mathfrak{d}_{F(v)/E}) \operatorname{disc}(E/\mathbb{Q})^2,
$
where $\mathfrak{d}_{F(v)/E}$ is the relative discriminant of the extension $F(v)/E$. This proves in particular that $\operatorname{disc}(E)^2$ divides $\operatorname{gcd}_{v \in \EPlaces} \operatorname{disc} F(v)$. Thus to finish the proof we need to show that for every prime $p$ we can find a place $v \in \EPlaces$ such that $N_{E/\mathbb{Q}}(\mathfrak{d}_{F(v)/E})$ is prime to $p$, or equivalently, that every place of $E$ of characteristic $p$ is unramified in $F(v)/E$.
Recall that $K'=K^{\text{conn}}$ is the minimal extension of $K$ over which all the endomorphisms of $A$ are defined, and that $K'$ is at most quadratic over $K$ (§\ref{sect_MT}). We now work with the abelian variety $A/K'$, so that the Galois representations have connected image and we can again apply theorem \ref{thm_Zarhin}.

Recall that we want to find a place $v$ in $\EPlaces$ such that $p$ is unramified in $F(v)/E$. We shall find $v$ as the restriction to $K$ of a suitable place $w$ of $K'$.
Let $(1,\omega)$ be a $\mathbb{Z}_p$-basis of $\mathcal{O}_E \otimes \mathbb{Z}_p$ and let $\omega^2-a\omega-b$ be the minimal polynomial of $\omega$. Notice that $\mathcal{O}_E \otimes \mathbb{Z}_p$ is not necessarily an integral domain, so the minimal polynomial of $\omega$ needs not be irreducible. 

Since $\mathcal{O}_E \otimes \mathbb{Q}_p$ acts on $V_p A=T_p A \otimes \mathbb{Q}_p$ in a way that is compatible with the action of Galois, we can identify $V_pA$ with $(\mathcal{O}_E \otimes \mathbb{Q}_p)^2$ as $(\mathcal{O}_E \otimes \mathbb{Q}_p)[G_{p^\infty}]$-modules, where the action of $G_{p^\infty} \subseteq \operatorname{GL}_2(\mathcal{O}_E \otimes \mathbb{Q}_p)$ on $(\mathcal{O}_E \otimes \mathbb{Q}_p)^2$ is the natural one.
 We can furthermore use the basis $(1,\omega)$ to identify $\mathbb{Q}_p^4$ with $(\mathcal{O}_E \otimes \mathbb{Q}_p)^{2}$ via $(x,y,z,w) \mapsto (x+y\omega,z+w\omega)$. We take
\[
f_p=
\begin{pmatrix}
 s & b & 0 & 0 \\
 1 & s+a & 0 & 0 \\
 0 & 0 & -a-s & b \\
 0 & 0 & 1 & -s \\
\end{pmatrix},
\]
where $s \in \mathbb{Q}_p$ is a parameter we will now choose: we claim that for all but finitely many $s \in \mathbb{Q}_p$ one has that $f_p$ is an element of $\mathcal{G}_p^0(\mathbb{Q}_p)$ whose characteristic polynomial has no multiple roots and whose centralizer in $\operatorname{End}_{\mathbb{Z}_p}(T_p A)$ is isomorphic to $(\mathcal{O}_E \otimes \mathbb{Z}_p)^2$. We now check these statements.

\begin{enumerate}
\item Recall our description of $\mathcal{G}_p^0$ as $\{ x \in \operatorname{GL}_2(E \otimes \mathbb{Q}_p) : \det(x) \in \mathbb{Q}_p^\times \}$ and the identification of $\mathbb{Q}_p^4$ with $(\mathcal{O}_E \otimes \mathbb{Q}_p)^2$. Via this identification, the operator $f_p$ is nothing but the $(\mathcal{O}_E \otimes \mathbb{Q}_p)$-linear operator corresponding to the matrix
\[
\begin{pmatrix}
s+\omega & 0 \\ 0 & -s-a+\omega
\end{pmatrix} \in \operatorname{GL}_2(\mathcal{O}_E \otimes \mathbb{Q}_p),
\]
whose determinant is $(s+\omega)(-s-a+\omega)=-s^2-as+b$, which is an element of $\mathbb{Q}_p$ for all $s \in \mathbb{Q}_p$ and is nonzero for all but finitely many of them. For any such value of $s$, the fact that the determinant is in $\mathbb{Q}_p^\times$, together with a short computation to check that $f_p$ is indeed a symplectic similitude, shows that $f_p$ lies in $\mathcal{G}_p^0(\mathbb{Q}_p)$.
\item The centralizer of $f_p$ is easily determined by a direct computation; it is isomorphic to $(\mathcal{O}_E \otimes \mathbb{Z}_p)^2$ provided that $a+2s$ and $s^2+as-b$ are nonzero, which is again true for all but finitely many values of $s$.
\item The discriminant of the characteristic polynomial of $f_p$ is $16 \left(a^2+4 b\right)^2 (a+2 s)^4 \left(a s-b+s^2\right)^2$,
so for all but finitely many values of $s$ the characteristic polynomial of $f_p$ has no repeated roots (if $a^2+4b$ were equal to zero, then we would have $(\omega-a/2)^2=0$, which would imply $\omega=a/2$ since $\mathcal{O}_E \otimes \mathbb{Z}_p$ is a reduced ring)
\end{enumerate}
 
Fix any value of $s$ for which $f_p$ has all the properties above.
By theorem \ref{thm_Zarhin} the set of places $w$ of $K'$ such that $\operatorname{End}_{\mathbb{F}_w}((A_{K'})_w) \otimes \mathbb{Z}_p \cong (\mathcal{O}_E \otimes \mathbb{Z}_p)^2$ has positive density. By theorem \ref{thm:Achter} and proposition \ref{prop:DensityOrdinaryReduction}, removing a 0-density subset of these we can further assume that $A_w$ is ordinary and absolutely irreducible for every such $w$.
Thanks to corollary \ref{cor_ordImpliesEndDefBase} we know that for such each $w$ we have
$
\operatorname{End}_{\mathbb{F}_{v}}(A_v) = \operatorname{End}_{\mathbb{F}_w}((A_{K'})_w),
$
where we have denoted by $v$ the place of $K$ induced by $w$. Notice furthermore that $A_v$ is ordinary and absolutely irreducible, so $v$ belongs to $\EPlaces$. We have thus shown that for a positive-proportion subset of $\EPlaces$ we have $\operatorname{End}_{\mathbb{F}_{v}}(A_v) \otimes \mathbb{Z}_p \cong (\mathcal{O}_E \otimes \mathbb{Z}_p)^2$. Since $\operatorname{End}_{\mathbb{F}_{v}}(A_v)$ is an order in $F(v)$ and $\mathcal{O}_E \otimes \mathbb{Z}_p$ is integrally closed, this implies in particular $\mathcal{O}_{F(v)} \otimes \mathbb{Z}_p \cong (\mathcal{O}_E \otimes \mathbb{Z}_p)^2$.
By lemma \ref{lemma_ANT1} we conclude that for every such $v$ the prime $p$ is unramified in the extension $F(v)/E$, which is what we needed to prove.
\end{itemize}
\end{proof}

Theorem \ref{thm:DetermineF} suggests the following procedure for determining $E$. The method is only guaranteed to give the correct result when $B \to \infty$, but in practice it produces the correct answer very quickly.
\begin{algorithm}\label{algo:DiscBound}
Let $A/K$ be an abelian surface for which $\operatorname{End}_{\overline{K}}^0(A)$ is a field. Fix a bound $B$.
\begin{enumerate}
\item compute $f_v(x)$ for all places $v \in \EPlaces$ with $q_v \leq B$.
\item for each $v$, compute the discriminant $\Delta(v)$ of the field generated by a root of $f_v(x)$.
\item if $\Delta(v)$ takes on at least two different values for different places $v$, output ``$A$ does not admit potential complex multiplication".
\item compute $d(B)=\operatorname{gcd}_{\substack{q_v \leq B \\ v \in \Omega_K'}} \Delta(v)$. If this number is 1, output ``$\operatorname{End}_{\overline{K}}^0(A)=\mathbb{Q}$ (hence $\operatorname{End}_{\overline{K}}(A)=\mathbb{Z}$)" and terminate the algorithm.
\item output ``the discriminant of $\operatorname{End}_{\overline{K}}^0(A)$ divides $d(B)$". If in step (3) we have proved that $A$ does not admit potential CM, output ``the square of the discriminant of $\operatorname{End}_{\overline{K}}^0(A)$ divides $d(B)$".
\end{enumerate}
\end{algorithm}

\begin{remark}\label{rem:FasterDiscBound}
One can in fact quit when $d(B)\leq 5^2-1$: indeed, the minimal discriminant of a quartic CM field is $125=\operatorname{disc}(\mathbb{Q}(\zeta_5))$, and the minimal discriminant of a real quadratic field is $5=\operatorname{disc}(\mathbb{Q}(\sqrt{5}))$.
In practice, this gives a significant speed improvement, because it is often the case that removing the last spurious factors of 2 from $d(B)$ requires taking $B$ much larger than what is needed to get rid of any other unwanted factor.
\end{remark}

\begin{remark}\label{rmk:EndKField}
There is an obvious variant of this algorithm that gives a bound on $\operatorname{disc}(\operatorname{End}_K(A))$, namely we compute $\Delta(v)$ for all places of good reduction for which $f_v(x)$ is irreducible. Since $\operatorname{End}_K(A) \hookrightarrow \operatorname{End}_{\mathbb{F}_v}(A_v)$, this gives an upper bound on $\operatorname{disc} (\operatorname{End}_K(A))$. However, it should be noted that it is not in general true that
\[
\operatorname{disc}(\operatorname{End}_K(A))= \operatorname{gcd}_{\substack{v \in \Omega_K^{(A)}}} \operatorname{disc}(F(v)):
\]
see example \ref{ex:DiscDoesNotWorkOverK} below. It seems likely that, for any given prime $q$, one can construct an abelian surface $A$ over a number field $K$ such that $\operatorname{End}_K(A)=\mathbb{Z}$ but $q \mid \operatorname{gcd}_{\substack{v \in \Omega_K^{(A)}}} \operatorname{disc}(F(v))$; however, such examples should be extremely rare.
\end{remark}

The output of algorithm \ref{algo:DiscBound} is correct because of theorem \ref{thm:DetermineF} and its proof. Moreover, choosing some reasonable bound $B$ ($B=200$ proved to be a good choice in our tests), we expect that the algorithm will correctly identify whether $\operatorname{End}_{\overline{K}}^0(A)$ is $\mathbb{Q}$, a real quadratic field, or a CM field, and that the bound of step (5) will be sharp. Notice that if $A$ has potential CM, then the field of complex multiplication is $F(v)$, where $v$ is any place in $\EPlaces$. Also notice that a real quadratic field is uniquely determined by the square of its discriminant, so the output of step (5) gives a reasonable candidate for $\EndKBar$.

To see why we expect the algorithm to give the correct answer quickly, notice first that already computing $f_v(x)$ for \textit{one} place in $\Omega_K'$ gives a nontrivial upper bound $d$ on the discriminant of $E$. Furthermore, for every prime $p$ dividing $d/\operatorname{disc}(E)$ (respectively $d/\operatorname{disc}(E)^2$ if $E$ is a quadratic field) the density of places $v$ for which $\operatorname{gcd}(d,\Delta(v))<d$ is positive, say it's $\delta>0$. At least heuristically, testing all the $N(B)$ places of norm up to $B$ then has probability $1-(1-\delta)^{N(B)}$ of removing one factor of $p$ from our estimate $d(B)$: since the distance of this quantity from 1 decays exponentially with $N(B)$, we expect the algorithm to eliminate any spurious factors in $\delta(B)$ fairly quickly.

In any case, the output of this algorithm leaves only a finite list of possibilities for the field $\EndKBar$. We remark that at this stage we have no information on the integral structure of $\operatorname{End}_{\overline{K}}(A)$, but if $\EndKBar$ is a field the ring $\operatorname{End}_{\overline{K}}(A)$ is certainly contained in the ring of integers of $\EndKBar$. Thus, provided that $\EndKBar$ is a field, we have obtained a finite list of rings $R_1, \ldots, R_k$ such that $\operatorname{End}_{\overline{K}}(A)$ embeds in one of the $R_i$ in such a way that the embedding has finite cokernel.

\section{Certifying the existence of extra endomorphisms}\label{sect:Certify}
The methods of the previous sections yield practical algorithm to obtain an ``upper bound" for the ring $\operatorname{End}_K(J)$ (or $\operatorname{End}_{\overline{K}}(J)$), namely a ring $R$ (or a finite list of rings $R_1,\ldots,R_k$) in which $\operatorname{End}_K(J)$ or $\operatorname{End}_{\overline{K}}(J)$ is contained (with finite index).
 We now focus on the converse problem of proving that a given two-dimensional Jacobian has nontrivial endomorphisms, at least over the algebraic closure. We have already discussed the reducible case in section \ref{sect_Reducibility}, so here we just consider the absolutely simple case. 

\subsection{Real multiplication}\label{sect:RM}
Assume that $J_{\overline{K}}$ has real multiplication. Letting the algorithm of section \ref{sect_Field} run with a sufficiently large bound $B$ (which is usually very small in practice) produces a list of real quadratic fields $E_1,\ldots, E_k$ (in all of our tests we found $k=1$) and a proof of the following statement: either $\operatorname{End}_{\overline{K}}(J)=\mathbb{Z}$, or $\operatorname{End}_{\overline{K}}(J)$ is an order in one of the $E_i$. Thus it suffices to check, for every $E_i$, whether $\operatorname{End}_{\overline{K}}(J)$ is an order in $E_i$. Consider a fixed $E=E_i$ and let $d$ be its discriminant.

The stratification of $\mathcal{A}_2$ described in section \ref{sect:ModuliSpace} implies that $\operatorname{End}_{\overline{K}}(J)$ is an order in $E$ if and only if the point $x_J \in \mathcal{A}_2$ corresponding to $J$ belongs to $\bigcup_{n \geq 1} H_{n^2d}$, and furthermore we have $x_J \in H_{n^2d}$ if and only if $\operatorname{End}_{\overline{K}}(J)$ is the unique order of index $n$ in $\mathcal{O}_E$. Thus if we can show that $x_J$ belongs to a certain Humbert surface $H_{n^2d}$ we have \textit{proved} the equality $\operatorname{End}_{\overline{K}}(J)=\mathcal{O}_{n^2d}$. The good news is that computing equations of Humbert surfaces can be reduced to a problem in linear algebra (see \cite{Gruenewald}), and furthermore the computation of such surfaces is completely independent of the curve $C$ we are considering, hence (even though it is a computationally expensive problem) we can consider that we have precomputed a large number of Humbert surfaces, and we only compute new ones whenever the need arises. Notice that a Humbert surface is defined by a single homogeneous polynomial in the Igusa invariants $I_2,I_4,I_6,I_{10}$. Testing whether $x_J \in H_{n^2d}$ is then immediate, because the Igusa invariants are polynomial functions of the coefficients of the curve, and we then just need to evaluate a single polynomial to check whether it vanishes on $x_J$.

The bad news, from the point of view of concrete implementation, is that the computation of equations for Humbert surfaces in $\mathbb{P}(I_2,I_4,I_6,I_{10})$ quickly becomes unfeasible. However, the situations becomes significantly better if one works in a certain finite-degree cover, the \textit{level-2 Satake model} $\mathcal{S}_2$ (with coordinates $x_1,\ldots,x_6$): points in $\mathcal{S}_2$ parametrize abelian varieties together with a fixed level 2 structure, and one can define Humbert surfaces $\tilde{H}_n$ in $\mathcal{S}_2$. 
For any point in $\mathcal{A}_2$ there are $|S_6|=720$ points in $\mathcal{S}_2$ mapping to it, and $x_J \in H_n$ if and only if at least one of the 720 inverse images belongs to $\tilde{H}_n$. From a computational point we then have the following procedure to check whether the point $x_J$ in moduli space corresponding to a certain Jacobian $J$ lies on $H_n$. First, one can compute the ``Satake coordinates" $s_i=\sum_{j=1}^6 x_j^i$ of $x_J$: these are just polynomial functions of the Igusa invariants $I_2,\ldots,I_{10}$ (see \cite[§3.4]{Gruenewald} for more details on level 2 Satake models). Then, using the classical relations between power sums and symmetric functions, from these one finds a degree-six polynomial whose six roots are the coordinates $x_1,\ldots,x_6$. Considering all permutations of these six coordinates gives the 720 points $x_{J,1}, \ldots, x_{J,720}$ in $\mathcal{S}_2$ that map to $x_J$, and for each of them we test whether $x_{J,i}$ belongs to $\tilde{H}_n$. For a single curve and a single Humbert surface (whose equation in $\mathcal{S}_2$ is known), this takes only a fraction of a second.
Furthermore, Gruenewald \cite{Gruenewald} has computed equations for all Humbert surfaces of discriminants up to 40, and for even discriminants up to 52, in the level 2 Satake model. This has been enough for all of our tests.

It is also worth mentioning that \textit{birational parametrizations} of certain degree-2 covers of the Humbert surfaces of primitive discriminant up to 100, and of square discriminant up to 121, are available thanks to the work of Elkies and Kumar \cite{MR3298543} \cite{MR3427148}. In principle, via Gr\"obner bases techniques these can also be used to test whether a point belongs to a Humbert surface. Of course this is incomparably slower than using the equations of $H_n$ when these are available, but it allows us to push the computation a bit further if necessary.

Finally, we remark that (as already hinted at in section \ref{sect_Reducibility}), one can also use equations of Humbert surfaces to prove that $J_{\overline{K}}$ is $(n,n)$-isogenous to a product of elliptic curves, simply by testing whether $x_J$ belongs to $H_{n^2}$. This can be used to quickly confirm geometric reducibility without going through the precedure of section \ref{sect_Reducibility}. Furthermore, this can also be used to confirm that the isogeny found by the method of §\ref{sect_Reducibility} has the minimal possible degree: indeed, $x_J \in H_{n^2}$ implies that $J$ is $(n,n)$ split, and is not $(m,m)$ split for any $m<n$.

\begin{remark}\label{rmk:Correspondences}
Recent work (see for example \cite{MR3540944}) has been directed toward showing that certain Jacobians admit actions of quadratic rings $\mathcal{O}_D$ \textit{over number fields}, that is, that there are embeddings $\mathcal{O}_D \hookrightarrow \operatorname{End}_K(J)$. Such work rely on entirely different techniques, namely, exhibiting explicit correspondences on the Jacobian, found by numerical methods not unlike those of section \ref{sect_Reducibility}.
\end{remark}

\subsection{Quaternionic multiplication}\label{sect:QM}

The problem of certifying quaternionic multiplication does not seem to have been extensively studied in the literature. We show that it can be reduced to that of certifying a finite number of real multiplications.

The idea is that an abelian variety has QM by a certain order $R$ in a quaternion algebra $Q$ if and only if it has real multiplication by a finite number of real quadratic rings:
\begin{proposition}\label{prop:CertifyQM}
Let $R$ be an order in a quaternion algebra $Q$. There exist finitely many discriminants $\Delta_1, \ldots, \Delta_n$ with the following property: let $R'$ be any quaternionic order. If there are optimal embeddings $\mathcal{O}_{\Delta_i} \hookrightarrow R'$ for all $i=1,\ldots,n$, then $R' \cong R$.
\end{proposition}
\begin{proof}
Let $M_R:=\begin{pmatrix} a & b \\ b & c \end{pmatrix}$ be a matrix associated with the order $R$ as in section \ref{sect:QuatRings}. Up to a $\operatorname{GL}_2(\mathbb{Z})$-change of basis, for example given by a matrix of the form $\begin{pmatrix}
k+1 & 1 \\ k & 1
\end{pmatrix}$, we can (and will) assume that $a$ and $c$ are distinct.
 Since $a$ and $c$ are both primitively represented by $M_R$, this implies that there are optimal embeddings $\mathcal{O}_a \hookrightarrow R$ and $\mathcal{O}_c \hookrightarrow R$. Consider now the set $\tilde{\mathcal{M}}(a,c)$ consisting of all the integral, positive-definite matrices of the form
\[
M_{a,c,x}=\begin{pmatrix} a & x \\ x & c \end{pmatrix};
\]
the condition that $M$ be positive-definite implies $0 \leq |x| \leq \sqrt{ac}$, hence $\tilde{\mathcal{M}}(a,c)$ is a finite set. Let $\sim$ be the equivalence relation on $\tilde{\mathcal{M}}_{a,c}$ given by ``$M_1 \sim M_2$ if and only if they represent the same quadratic form", and let $\mathcal{M}_{a,c}$ be a set containing $M_R$ and consisting of one representative for each equivalence class in $\tilde{\mathcal{M}}_{a,c}/\sim$.

The classical theory of quadratic forms implies that no two elements of $\mathcal{M}(a,c)$ primitively represent the same set of integers, hence for every $M \in \mathcal{M}(a,c) \setminus \{ M_R\}$ we can find an integer $D=D(M)>0$ that is primitively represented by $M_R$ but not by $M$ (recall that $M$ is positive definite, hence the condition $D > 0$ is automatic). In particular, there is an optimal embedding of $\mathcal{O}_D$ in $R$.
We claim that we can take 
\[
\{\Delta_1,\ldots,\Delta_n\}=\{a,c\} \cup \bigcup_{M \in \mathcal{M}(a,c) \setminus \{M_R\}} \{D(M)\}.
\]
We have already seen that all the corresponding quadratic rings $\mathcal{O}_{\Delta_i}$ optimally embed in $R$, so it suffices to show the other implication.

Suppose $R'$ is a quaternionic order in which both $\mathcal{O}_a$ and $\mathcal{O}_c$ embed optimally. Then, up to $\operatorname{GL}_2(\mathbb{Z})$-equivalence and using the fact that $a \neq c$, its discriminant matrix $M_{R'}$ is of the form $\begin{pmatrix} a & x \\ x & c \end{pmatrix}$, hence in particular equivalent to an element of $\mathcal{M}(a,c)$. If by contradiction we had $M_{R'} \neq M_R$ (as quadratic forms) there would exist an integer $D=D(M_{R'}) \in \{\Delta_1,\ldots,\Delta_n\}$ that is not primitively represented by $M_{R'}$: but this implies (by theorem \ref{thm_RMImpliesEverything} (4)) that there is no optimal embedding $\mathcal{O}_{D} \hookrightarrow R'$, contradiction. 
\end{proof}

Modifying slightly the proof of the previous proposition we obtain the following result:

\begin{theorem}\label{thm:CertifyQM}
Let $R$ be a quaternionic order in a quaternion algebra over $\mathbb{Q}$. There exist two finite, and explicitly computable, sets of discriminants $\mathcal{P}_R=\{P_1,P_2\}$ and $\mathcal{N}_R=\{N_1,\ldots,N_n\}$ with the following property: let $A/K$ be a principally polarized abelian variety. Then $A$ admits an optimal action of $R$ defined over $K$ if and only if the following two conditions are met:
\begin{enumerate}
\item for $D=P_1,P_2$, there is an optimal action of $\mathcal{O}_D$ on $A$ defined over $K$;
\item for each $D \in \mathcal{N}_R$ there is no optimal action of $\mathcal{O}_D$ on $A$ defined over $\overline{K}$; equivalently the point in $\mathcal{A}_2$ corresponding to $A$ does not lie on the Humbert surface $H_D$.
\end{enumerate}
\end{theorem}
\begin{proof}
Reasoning as in the previous proposition, let $P_1=a,P_2=c$ be discriminants such that $\mathcal{O}_a, \mathcal{O}_c$ optimally embed in $R$. For each $M_i \in \mathcal{M}_{a,c} \setminus \{M_R\}$, let $N_i$ be an integer which is primitively represented by $M_i$ but not by $M_R$ (any such integer is automatically a discriminant, as remarked in section \ref{sect:QuatRings}). Now conditions (1) and (2) clearly hold if $R$ embeds optimally in $\operatorname{End}_K(A)$; conversely, suppose that (1) and (2) hold. Since we have two optimal embeddings $\mathcal{O}_a \hookrightarrow \operatorname{End}_K(A), \mathcal{O}_c \hookrightarrow \operatorname{End}_K(A)$ we deduce that $\operatorname{End}_K(A)$ contains a quaternion ring. Up to a change of basis the discriminant matrix of this ring is an element of $\mathcal{M}_{a,c}$, and condition (2) guarantees that it must in fact be $M_R$.
\end{proof}

\begin{remark}
In section \ref{sect:RM} we have discussed using explicit equations of Humbert surfaces to determine whether or not there is an optimal action of $\mathcal{O}_D$ on $A_{\overline{K}}$, which is precisely what we need to verify condition (1) in the previous theorem. The work mentioned in remark \ref{rmk:Correspondences} is also relevant here: the previous theorem implies that any technique useful to certify that a genus-2 Jacobian admits an action of certain quadratic rings (over a certain field $K'$) can also be used to demonstrate the existence of quaternionic multiplication on that same Jacobian (again over $K'$). Notice however that we need to know that the action is \textit{optimal}; suppose that we only know that there is an action of $\mathcal{O}_D$ on $A$ defined over $K$, but not necessarily that it is optimal. A useful criterion is the following: if for every positive integer $n>1$, $n^2 \mid D$, there is no action of $\mathcal{O}_{D/n^2}$ on $A$, then the given action of $\mathcal{O}_D$ is optimal. Hence (if we know that there is an action of $\mathcal{O}_D$ on $A/K$) to prove optimality it suffices to check that the moduli point corresponding to $A$ does not lie on $\bigcup_{\substack{n^2>1\\ n^2 \mid D}}H_{D/n^2}$. Obviously, if $D$ is squarefree the action is always optimal.
\end{remark}

\begin{remark}
Theorem \ref{thm:CertifyQM} is more practical than proposition \ref{prop:CertifyQM} from a computational point of view, essentially because it reduces the number of optimal embedding one needs to find from the \textit{a priori} arbitrary number $n$ to just 2.
The price is pay is that we also need to prove that $A$ does \textit{not} admit optimal actions of certain quadratic rings. This, however, should be considered easier than the problem of showing that a certain quadratic ring does act optimally on $A$: one reason for this is that \textit{not} admitting an optimal action of a ring $R$ is an open condition in the moduli space, and is therefore amenable to approximate computations. Moreover, one can also try to obtain \textit{negative} information on $\operatorname{End}_{\overline{K}}(A)$ by looking at the reductions of $A$ (while it is hard to prove something \textit{positive} about $\operatorname{End}_{\overline{K}}(A)$ by just considering reductions).
\end{remark}

\subsection{Complex Multiplication}\label{sect:CM}
All the genus 2 curves defined over the rational field whose Jacobian is geometrically simple and admits (potential) CM by the maximal order of a quartic field have been listed in \cite{MR1609658} (that the list is complete has been proven in \cite{MR1807666}). The list of such curves is finite, and in fact pretty short, so we can test whether we are in either of these cases simply by matching the absolute Igusa invariants of our curve against this list. 
We shall not discuss this case further, because an algorithmic approach to proving that a given Jacobian has complex multiplication has already been discussed in \cite{MR1648415}. We do remark, however, that in our tests we did not find any genus-2 curves with complex multiplication by a non-maximal order in a quartic field.

\subsection{Real multiplication over number fields: a remark}\label{sect:PotentialRM}
It is worth out mentioning that in the special case of \textit{potential} real multiplication our methods often allow us to determine the structure of $\operatorname{End}_{\overline{K}}(A)$ even as a Galois module.
Indeed, let $A/K$ be an abelian surface. Suppose that we have proved that $\EndKBar$ is a field, and -- using the method of section \ref{sect:RM} -- we have then concluded that $\operatorname{End}_{\overline{K}}(A)$ is a certain order $R$ in a real quadratic field. Suppose furthermore that using Algorithm \ref{algo:DiscBound} (in the version that computes an upper bound for $\operatorname{End}_K(A)$, i.e. remark \ref{rmk:EndKField}) we have been able to prove $\operatorname{End}_K(A)=\mathbb{Z}$, so that the real multiplication is only potential. In order to determine $\operatorname{End}_{\overline{K}}(A)$ as a Galois module, it now suffices to compute the action of $\abGal{K}$ on $R$.

By section \ref{sect_MT} we know that there exists an \textit{at most quadratic} extension $K'$ of $K$ over which all the endomorphisms of $A$ are defined, and since by assumption $K' \neq K$ we have $[K':K]=2$. The action of $\abGal{K}$ on $R$ plainly factors through $\operatorname{Gal}(K'/K)$, and the nontrivial element of $\operatorname{Gal}(K'/K)$ necessarily acts on $R$ as its unique nontrivial involution that fixes $\mathbb{Z}$. Thus the only ingredient left to be determined is the field $K'$.

It is known that the places ramified in $K'/K$ are a subset of the places where $A$ does not have semistable reduction (this was noticed by Ribet, see \cite[Page 262]{MR1154704}), so given $A$ we can compute a finite list of candidate fields $K'$, namely those quadratic extensions of $K$ unramified outside $\{v \in \Omega_K : A_v \text{ is not semistable} \}$. We remark that if $A$ is the Jacobian of the genus 2 hyperelliptic curve $y^2=f(x)$ with $f(x) \in \mathcal{O}[x]$, where $\mathcal{O}$ is the ring of integers of $K$, then a sufficient condition for $A$ to have semistable reduction at a place $v \in \Omega_K$ of characteristic not 2 is that the $v$-adic valuation of $\operatorname{disc} f(x)$ is at most 1 (see for example Kowalski's appendix in \cite{MR2820155}).

Let $K_1', \ldots, K_m'$ be the finite list of candidate fields $K'$. Notice that if $K'_i \neq K'$ then $\operatorname{End}_{K'_i}(A)=\mathbb{Z}$.
For each field $K_i'$ we can then use algorithm \ref{algo:DiscBound} (in the version for $\operatorname{End}_K(A)$, see remark \ref{rmk:EndKField}) to compute an upper bound on $\operatorname{End}_{K_i'}(A)$. Choosing a bound $B$ large enough we can hope to show that for $m-1$ candidate fields we have $\operatorname{End}_{K_i'}(A)=\mathbb{Z}$, thus proving $K' \neq K_i'$. The only remaining candidate field is now \textit{provably} the minimal field over which the endomorphisms of $A$ are defined.

A small variant of this idea can be turned into a deterministic algorithm as follows. Notice that in the case of potential real multiplication the set of places $v$ for which $A_v$ is absolutely irreducible and ordinary has density 1 (theorem \ref{thm:Achter} and proposition \ref{prop:DensityOrdinaryReduction}). For every such $v$ we have a corresponding quartic CM field $F(v)$, and we denote by $F_0(v)$ the unique real quadratic subfield of $F(v)$. Notice that if $\operatorname{End}^0_{K_i}(A)$ is a real quadratic field then for every such $v$ we have $F_0(v)=\operatorname{End}^0_{K_i}(A)$ (see e.g.~the proof of theorem \ref{thm:DetermineF}), so, in order to show that $A$ does not admit real multiplication over $K_i$, it suffices to prove that for some place $v$ of $K_i$ (of good ordinary reduction for $A$, with $A_v$ absolutely simple) the fields $F_0(v)$ and $\EndKBar$ are not isomorphic. We now show that such a place $v$ always exists, which -- combined with our previous remarks -- leads to a (fairly efficient) algorithm to determine the smallest field of definition of the real multiplication of $A$, see algorithm \ref{algo:FieldDefinitionRM} below.
\begin{proposition}\label{prop:ComputeFieldDefinitionRealMultiplication}
Let $A/K$ be an abelian surface such that $E=\EndKBar$ is a real quadratic field and $\operatorname{End}_K(A)=\mathbb{Z}$. The places $v \in \Omega_K^{(A)}$ such that $A_v$ is ordinary and absolutely simple, and $F_0(v)$ is not isomorphic to $E$, have positive density.
\end{proposition}
\begin{proof}
There is a bound $\ell_0$ such that for every prime $\ell>\ell_0$ the group $G_{\ell^\infty}$ contains $G_{\ell^\infty}^0 := \{x \in \operatorname{GL}_2(\mathcal{O}_E \otimes \mathbb{Z}_\ell)\bigm\vert \det(x) \in \mathbb{Z}_\ell^\times \}$ as an index-2 subgroup (see §\ref{sect_MT}). Pick a prime $\ell > \max\{3,\ell_0\}$ that is nonsplit in $E$ and fix a basis $1,\omega$ of $E \otimes \mathbb{Z}_\ell$ over $\mathbb{Z}_\ell$ such that $\omega^2=d \in \mathbb{Z}_\ell^\times$ (this choice of $\omega$ can be made since $\ell \neq 2$). The $\mathbb{Z}_\ell$-module $T_\ell(A)$ is a free $(\mathcal{O}_E \otimes \mathbb{Z}_\ell)$-module of rank 2; fix a $(\mathcal{O}_E \otimes \mathbb{Z}_\ell)$-basis $x,y$ and use $x,y,\omega x,\omega y$ as $\mathbb{Z}_\ell$-basis of $T_\ell(A)$. With respect to these coordinates, the matrices in $G_{\ell^\infty}^0$ can be represented as
\[
\begin{pmatrix}
a_{11} & a_{12} & d b_{11} & d b_{12} \\
a_{21} & a_{22} & d b_{21} & d b_{22} \\
b_{11} & b_{12} & a_{11}   & a_{12} \\
b_{21} & b_{22} & a_{21}   & a_{22}
\end{pmatrix}, \quad a_{11}b_{22}-a_{12}b_{21}-a_{21} b_{12}+a_{22}b_{11}=0,
\]
and those in $G_{\ell^\infty} \setminus G_{\ell^\infty}^0$ as 
\[
\begin{pmatrix}
a_{11} & a_{12} & d b_{11} & d b_{12} \\
a_{21} & a_{22} & d b_{21} & d b_{22} \\
-b_{11} & -b_{12} & -a_{11}   & -a_{12} \\
-b_{21} & -b_{22} & -a_{21}   & -a_{22}
\end{pmatrix}, \quad a_{11}b_{22}-a_{12}b_{21}-a_{21} b_{12}+a_{22}b_{11}=0,
\]
where the coefficients $a_{ij}, b_{ij}$ are in $\mathbb{Z}_\ell$, the determinant of the $4 \times 4$ matrices is an $\ell$-adic unit, and the bilinear equation corresponds to the condition $\det(x) \in \mathbb{Z}_\ell^\times$.
Now we observe that $G_{\ell^\infty} \setminus G_{\ell^\infty}^0$ contains operators whose characteristic polynomial splits completely in $\mathbb{F}_\ell[x]$ and has no repeated roots modulo $\ell$ (for example one can set $b_{ij}=0$ for $i,j =1,2$, $a_{12}=a_{21}=0$, $a_{11}=1$, $a_{22}=2$; the roots of the characteristic polynomial are then distinct modulo $\ell>3$). Let $C$ be the conjugacy set in $G_\ell$ given by those operators whose characteristic polynomial is separable and split over $\mathbb{F}_\ell$; by the above, it is non-empty.
By Chebotarev's density theorem, there is a set $T$ of places $v \in \Omega_K^{(A)}$ of positive density such that $\rho_\ell(\operatorname{Fr}_v)$ lies in $C$ for every $v \in T$. Since places $v$ such that $A_v$ is absolutely simple have density 1 (theorem \ref{thm:Achter}), and the same is true for places of ordinary reduction (proposition \ref{prop:DensityOrdinaryReduction}), the set
\[
T' := \{ v \in T : A_v \text{ is absolutely simple and ordinary} \}
\]
has positive density as well.
For every $v \in T'$ we have that 
$f_v(x)$ is split and separable modulo $\ell$: this implies in particular that $\ell$ is completely split in $F(v)$, hence also in $F_0(v)$. Since $\ell$ was chosen so as \textit{not} to split in $E$ this proves that $F_0(v)$ and $E$ are not isomorphic.
\end{proof}

We have thus obtained the following algorithm (which terminates thanks to proposition \ref{prop:ComputeFieldDefinitionRealMultiplication}):
\begin{algorithm}\label{algo:FieldDefinitionRM}
Let $A/K$ be an abelian surface such that $\EndKBar$ is a real quadratic field $E$ (assumed to be known) but $\operatorname{End}_K(A)=\mathbb{Z}$.
\begin{itemize}
\item Compute a finite list $L=\{K_i\}$ of quadratic extensions of $K$ such that the equality $\EndKBar=\operatorname{End}_{K_i}(A)$ holds for at least (hence precisely) one field $K_i$.
\item Loop over an integer $B$. For every field $K_i$ in $L$:
\begin{enumerate}
\item compute $f_v(x)$ for $v \in \Omega_{K_i}^{(A)}$ of norm at most $B$.
\item for every $v$ as above test whether $f_v^{[12]}(x)$ is irreducible and $A_v$ is ordinary. If these conditions hold:
\begin{enumerate}
\item compute the real subfield $F_0(v)$ of $F(v)=\frac{\mathbb{Q}[x]}{(f_v(x))}$.
\item if $F_0(v)$ and $E$ are not isomorphic, remove $K_i$ from $L$.
\end{enumerate}
\end{enumerate}
\item If $|L|=1$, the only field remaining in $L$ is the minimal extension of $K$ over which the real multiplication of $A$ is defined. If $|L|>1$, increase $B$.
\end{itemize}
\end{algorithm}

\begin{remark}\label{rmk:Congruences}
Suppose for simplicity that $K=\mathbb{Q}$ (but the same remark applies, mutatis mutandis, to arbitrary number fields). Running the algorithm of remark \ref{rmk:EndKField} over the quadratic field $K_i'$ is essentially equivalent to running it over $\mathbb{Q}$, but restricting to primes that split in $K_i'$: indeed, it is well-known that the set of places of degree 1 has full density, and moreover, nothing would be gained by looking at the inert places $w$, because since $[\mathbb{F}_w:\mathbb{F}_{p_w}]=2=[K':K]$ we know that the real multiplication is defined over $\mathbb{F}_w$. 

Also notice that for the primes $v$ of $K_i'$ of degree one the residue field is just the prime field, so that we might as well compute $f_{p_v}(x)$ for the prime $p_v$ lying under $v$. Finally, the primes $p_v$ that we need to consider are easy to determine, because they are given by congruence conditions: in other words, we can carry out our proposed test using only arithmetic in $\mathbb{F}_p$ for those primes $p$ that lie in certain explicit arithmetic progressions.
\end{remark}

\begin{example}\label{ex:RMOverGroundField}
Consider the curve $C: y^2 = x^5 - x^4 - x^3 + x^2 + x - 1$. Its conductor is $2^{12} 3^2$, so $A:=\operatorname{Jac}(C)$ has good reduction away from 2 and 3. Algorithm \ref{algo:Irreducibility} with $B=7$ proves that $A$ is geometrically irreducible without potential QM, so its absolute endomorphism ring is an order in a number field.
Algorithm \ref{algo:DiscBound} with $B=67$ (and even $B=23$, if we use Remark \ref{rem:FasterDiscBound}) shows that $\operatorname{End}_\mathbb{Q}(A)=\mathbb{Z}$. The same algorithm (in the version for $\EndKBar$) proves that $A$ does not have potential CM, and that $(\operatorname{disc} \EndKBar)^2$ divides 64, so we suspect that $A$ has potential real multiplication by an order in $\mathbb{Q}(\sqrt{2})$. Using the equation of the Humbert surface $H_8$, we determine that $\operatorname{End}_{\overline{K}}(A)=\mathbb{Z}[\sqrt{2}]$. The minimal field of definition of the real multiplication is a quadratic extension of $\mathbb{Q}$ unramified away from $2$ and $3$, so it is one of the following fields: $\mathbb{Q}(i), \mathbb{Q}(\sqrt{\pm 2}), \mathbb{Q}(\sqrt{\pm 3}), \mathbb{Q}(\sqrt{\pm 6})$. Finally, algorithm \ref{algo:DiscBound} with $B=61$ (and considering only primes in the correct congruence classes) proves that $\operatorname{End}_{F}A=\mathbb{Z}$ for $F=\mathbb{Q}(i), \mathbb{Q}(\sqrt{-2}), \mathbb{Q}(\sqrt{\pm 3}), \mathbb{Q}(\sqrt{\pm 6})$ -- notice that, as shown by the following example, for this computation we really need the improvement described in Remark \ref{rem:FasterDiscBound}.
Putting everything together, we have proved that the minimal field of definition of the real multiplication of $A$ is $\mathbb{Q}(\sqrt{2})$.
\end{example}

\begin{example}\label{ex:DiscDoesNotWorkOverK}
The previous example displays the behaviour described in Remark \ref{rmk:EndKField}, namely, even though $\operatorname{disc} \operatorname{End}_K^0(A)=\operatorname{disc} \mathbb{Q}=1$ for $K=\mathbb{Q}(\sqrt{-2})$, one has
\[
\operatorname{gcd}_{v \in \Omega_K^{(A)}} \operatorname{disc}(F(v)) \neq 1.
\]
We prove this claim. Notice first that as pointed out in remark \ref{rmk:Congruences} we only need consider the places $v$ of $K$ of degree 1; moreover, $A$ has bad reduction at the unique place of $\mathbb{Q}(\sqrt{-2})$ lying over 2, so we only need to consider primes that split completely in $K$. We denote by $v$ a degree-1 place of $K$ and by $p$ the rational prime below it; the characteristic polynomial of the Frobenius at $v$ is thus the same as the characteristic polynomial of the Frobenius at $p$ acting on $\operatorname{Jac}(C/\mathbb{Q})$.
A rational prime $p$ splits completely in $K$ if and only if it is congruent to $1,3 \pmod 8$. Those that are congruent to $1 \pmod 8$ also split in $\mathbb{Q}(\sqrt{2})$, so for such $v$ one has $\mathbb{Z}[\sqrt{2}] \hookrightarrow F(v)$ since the real multiplication is defined over $\mathbb{Q}(\sqrt{2})$, hence over $\mathbb{F}_v$. This shows that $2 \mid \operatorname{disc} F(v)$. We now consider primes $p$ congruent to $3$ modulo 8. These do not split in $\mathbb{Q}(\sqrt{2})$, so they must map to the nontrivial class in $\mathcal{G}_\ell/\mathcal{G}_\ell^0$: indeed, by lemma \ref{lemma:ExpDivides12} one has that the kernel of the natural map $\abGal{\mathbb{Q}} \to \mathcal{G}_\ell/\mathcal{G}_\ell^0$ is $\operatorname{Gal}\left( \overline{\mathbb{Q}} / \mathbb{Q}(\sqrt{2}) \right)$. Fix (for simplicity) $\ell=2$. Arguing as in the proof of Proposition \ref{prop:ComputeFieldDefinitionRealMultiplication}, we can find a basis of $T_2(A)$ such that 
the elements in the nontrivial coset $\mathcal{G}_2(\mathbb{Q}_2) \setminus \mathcal{G}_2^0(\mathbb{Q}_2)$ can be written as
\[
\begin{pmatrix}
 a_1 & b_1 & 2 a_2 & 2 b_2 \\
 c_1 & d_1 & 2 c_2 & 2 d_2\\
-a_2 & -b_2 & -a_1 & -b_1 \\
 -c_2 & -d_2 & -c_1 & -d_1 \\
\end{pmatrix}, \quad  \begin{array}{c}a_1,a_2,b_1,b_2,c_1,c_2,d_1,d_2 \in \mathbb{Q}_2 \\ \text{ with } a_2 d_1 + 
a_1 d_2=b_2 c_1 + b_1 c_2\end{array}.
\]
We conclude that the trace of any operator lying in $\mathcal{G}_2(\mathbb{Q}_2) \setminus \mathcal{G}_2^0(\mathbb{Q}_2)$ vanishes; in particular, the trace of the Frobenius at $v$ -- which can be computed as $\operatorname{tr} \rho_{2^\infty}(\operatorname{Fr}_v)$ -- is equal to zero
(the same conclusion can also be obtained from \cite[§3.6]{MR2982436}). It follows that the characteristic polynomial of the Frobenius at $v$ is of the form $x^4+ax^2+p^2$ for some integer $a$, and we now show that $a$ is even. One has 
\[
x^5-x^4-x^3+x^2+x-1=(x-1)(x^4-x^2+1),
\]
so that a basis of the 2-torsion of $A$ is given by the divisors
\[
P_1=(\zeta_{12},0)-\infty,  \quad P_2=(\zeta_{12}^5,0)-\infty, \quad P_3=(\zeta_{12}^7,0)-\infty, \quad P_4=(\zeta_{12}^{11},0)-\infty,
\]
where $\zeta_{12}$ is a primitive 12-th root of unity. This gives an explicit description of the action of $\abGal{K}$ on $A[2]$, from which one sees that every characteristic polynomial of Frobenius is congruent to $t^4+1$ modulo 2.
We conclude that for $p \equiv 3 \pmod 8$ the characteristic polynomial of $\operatorname{Fr}_v$ is of the form $x^4+2bx^2+p^2$ for some integer $b$, and this implies that 2 ramifies in $F(v)$. Thus $2 \mid \operatorname{gcd}_{v \in \Omega_K^{(A)}} \operatorname{disc}(F(v))$; in fact, with some more effort, one can show $\operatorname{gcd}_{v \in \Omega_K^{(A)}} \operatorname{disc}(F(v))=16$.
Finally, notice that we were somewhat lucky with this example, in the sense that $\operatorname{gcd}_{v \in \Omega_K^{(A)}} \operatorname{disc}(F(v))$ is strictly less than $5^2$, so that we can apply remark \ref{rem:FasterDiscBound}. If $\operatorname{gcd}_{v \in \Omega_K^{(A)}} \operatorname{disc}(F(v))$ had exceeded $5^2-1$ (as might happen with examples with larger coefficients), we would have had to resort to algorithm \ref{algo:FieldDefinitionRM} to conclude that $\operatorname{End}_K(A)=\mathbb{Z}$.
\end{example}

\section{Some numerical tests}\label{sect:NumericalResults} 

\subsection{Curves with small coefficients}

To test how the algorithm performs on a ``random" genus-2 curve with small coefficients, we have applied it to all models of the form $y^2=x^5+a_4x^4+a_3x^3+a_2x^2+a_1x+a_0$, where the $a_i$ are integers bounded by 10 in absolute value and $a_4 \geq 0$ (this restriction is imposed so as to avoid double-counting curves that only differ by an obvious hyperelliptic twist). We give some statistics.

\begin{itemize}
\item Out of the $2139291$ models thus tested, $7239$ are singular and have been excluded from the computation.
\item For $2130158$ of the remaining $2132052$ models, algorithm \ref{algo:Irreducibility} shows that the corresponding Jacobian is geometrically irreducible without potential QM. In all cases choosing $B=59$ has been sufficient for this algorithm (that is, we have never computed $f_v(x)$ for a prime $v$ of $\mathbb{Q}$ of norm 61 or more). 
\item Using equations of Humbert surfaces, we have checked that $1885$ of the remaining models have a Jacobian which is not geometrically irreducible, and in each case have determined the minimal $n$ such that $J$ is geometrically $(n,n)$-isogenous to a product of elliptic curves.
At this point, for 9 models we do not know yet whether $J$ is geometrically irreducible or not (hence we strongly suspect them to admit potential quaternionic multiplication).
\item Using the method of theorem \ref{thm:CertifyQM}, the 9 remaining models have been proven to have a Jacobian with potential QM. Thus in all cases we have \textit{proved} whether or not $\operatorname{Jac}(C)$ is geometrically irreducible.
\item Combining the previous results with the method of sections \ref{sect_Field} and \ref{sect:RM} we have computed the geometric endomorphism ring of all the curves considered. The results are summarized in the following table:

\begin{center}
\begin{tabular}{c|cc}
& $\operatorname{End}_{\overline{\mathbb{Q}}}$ & Number of models \\ \hline
Trivial & $\mathbb{Z}$ & 2129918 \\ \hline
CM & $\mathbb{Z}[\zeta_5]$ & 41 \\ 
 & $\mathbb{Z}[\sqrt{-2+\sqrt{2}}]$ & 1 \\ 
\hline
RM & $\mathbb{Z}\left[ \frac{1+\sqrt{5}}{2} \right]$ & 84 \\
& $\mathbb{Z}\left[ \sqrt{2} \right]$ & 95 \\
& $\mathbb{Z}\left[ \sqrt{3} \right]$ & 7 \\
& $\mathbb{Z}\left[ \sqrt{5} \right]$ & 2 \\
& $\mathbb{Z}\left[ \sqrt{6} \right]$ & 2 \\
& $\mathbb{Z}\left[ \sqrt{13} \right]$ & 2 \\
& $\mathbb{Z}\left[ \sqrt{17} \right]$ & 6 \\ 	\hline
Decomposable & $(2,2)$-decomposable & 1810 \\
& $(3,3)$-decomposable & 66 \\
& $(4,4)$-decomposable & 6 \\
& $(5,5)$-decomposable & 3 \\ \hline
QM & Maximal order of $Q_6$ & 3 \\ 
& $\begin{array}{c}\text{Order of index 2}\\ \text{ in a maximal order of }Q_6\end{array}$ & 6 \\ \hline
Total & &  2132052 \\
\end{tabular}
\end{center}
\end{itemize}

In the QM case, $Q_6$ denotes the unique $\mathbb{Q}$-quaternion algebra ramified at 2 and 3. The non-maximal order appearing in the next to last row is $\mathbb{Z} \oplus \mathbb{Z} \alpha \oplus \mathbb{Z} \beta \oplus \mathbb{Z} \alpha\beta$, where $\alpha^2=\beta^2=2$, $\beta\alpha+\alpha\beta=-2$ (to see that this is indeed an order in $Q_6$ it may be helpful to notice that $(2\alpha+\beta)^2=6$).

\subsection{The genus-2 curves in the LMFDB}\label{sect:LMFDB}

We have run our algorithm on all the genus-2 curves listed in the LMFDB whose Jacobians, according to the data computed in \cite{MR3540958}, have nontrivial geometric endomorphism ring. In all cases our results matched those in the LMFDB: for all the geometrically irreducible cases we have confirmed the structure of $\operatorname{End}_{\overline{\mathbb{Q}}}(\operatorname{Jac}(C))$, and for all the geometrically reducible cases we have checked that $\operatorname{Jac}(C)$ maps to the two correct elliptic curves with an isogeny of the correct degree.
We describe in detail the computations leading to the determination of $\operatorname{End}_{\overline{\mathbb{Q}}}(\operatorname{Jac}(C))$ for two interesting examples.

\subsubsection{Non-simple Jacobian}
We carry out the procedure of section \ref{sect_Reducibility} on the Jacobian $J$ of the 
curve labelled 20412.b.734832.1, namely
\[
C : y^2+(x^2+x)y=x^6+3x^5+2x^4+7x^3+11x^2+14;
\]
according to the LMFDB, this curve admits a map towards the elliptic curve 54.a2 (Cremona label 54a3), that is, the curve with minimal Weierstrass equation
\[
E_{54.a2} : z^2+wz= w^3-w^2 -3w+3. 
\]
The curve $C$ is unique in its being the only one currently indexed by the LMFDB that admits a degree-7 map to an elliptic curve, and no maps of lower degree. We shall prove that there is indeed a degree 7 map to an elliptic curve by finding it explicitly.
We find it easier to work with the purely hyperelliptic models
$
C: y^2=4x^6 + 12x^5 + 9x^4 + 30x^3 + 45x^2 + 56
$
and
$
E_{54.a2} : z^2 = w^3 - 51w + 142.
$

From now on, all of our findings are purely numerical; we shall state them as facts, but the reader should keep in mind that we will only know that these computations are rigourously correct once we have found an explicit map $C \to E$.
MAGMA's intrinsic AnalyticEndomorphisms reports that $\operatorname{End}(J)$ is the $\mathbb{Z}$-subalgebra of $M_4(\mathbb{Z})$ generated (as an algebra) by the identity matrix and by 
$
M_1:=\begin{pmatrix}
  1 &  1 &  0 & -1 \\
 -2 &  0 &  1 &  0 \\
  0 & 14 &  1 & -2 \\
  -14 &  0 &  1 & 0
\end{pmatrix}.
$
We find that $M_1$ has two integral eigenvalues, namely 4 and $-3$, both of multiplicity 2. We then consider the rank-2 matrix
$
M:=M_1-4 \operatorname{Id} = \begin{pmatrix}
 -3 &  1 &  0 & -1 \\
 -2 & -4 &  1 &  0 \\
  0 & 14 & -3 & -2 \\
-14 &  0 &  1 & -4
\end{pmatrix}.
$
We know the corresponding $A$-matrix only numerically, but rational reconstruction gives $A=\begin{pmatrix}
-7 & -1 \\ 0 & 0
\end{pmatrix}$. The column space of $M$ is generated by its third and fourth column.

We thus obtain a numerical approximation to $\Pi M$, and the $\mathbb{Z}$-span of the columns of $\Pi M$ is generated by its third and fourth column (call them $C_3, C_4$). Let $\Lambda_{E,2}$ be the lattice generated in $\mathbb{C}^2$ by $C_3, C_4$. As projection $\pi : \mathbb{C}^2 \to \mathbb{C}$ we choose $(z_1,z_2) \mapsto z_1$, which sends $C_3, C_4$ to two $\mathbb{Z}$-linearly independent complex numbers $\omega_1,\omega_2$. Denote by $\Lambda_E$ the lattice generated in $\mathbb{C}$ by the images of $C_3, C_4$ via $\pi$. We compute numerically the $j$-invariant of $\mathbb{C}/\Lambda_E$, and (using rational reconstruction again) we find $j=-132651/2$. Computing $g_4(\Lambda_E), g_6(\Lambda_E)$ we obtain the equation of $E=E_{54.a2}$. In the notation of section \ref{sect_Reducibility}, the kernel of $\psi:\mathbb{C} \to \mathbb{C}/\Lambda_E$ contains $\tilde{\Lambda}$ with index 7 (this is easy to see from the matrix expression of $A$). We then look for a degree-7 map $C \to E$; proceeding as in section \ref{sect_Reducibility} we find a nontrivial map $(x,y) \mapsto (w(x),z(x,y))$ from $C$ to $E_{54.a2}$, an explicit expression for which is given by
\[
\displaystyle 
\begin{pmatrix}
w \\ z
\end{pmatrix}
 = 
 \begin{pmatrix} \displaystyle  \frac{1-78 x+72 x^2-102 x^3+36 x^4-12 x^5-7 x^6+6 x^7}{(1-x+x^2)^2 (7+3 x^2+2 x^3)} \\ \phantom{space} \\ \displaystyle  \frac{4 (-19+18 x-12 x^2+13 x^3+9 x^4+12 x^5+4 x^6+9 x^7+2 x^9)}{(1-x+x^2)^3 (7+3 x^2+2 x^3)^2}y\end{pmatrix},
\]
and which has degree 7 as predicted. Notice that checking that this map gives a covering $C \to E_{54.a2}$ amounts to some trivial algebra, so even if we cannot rigorously justify any of the previous computations the end result is provably correct.

We can also carry out the same computation for another elliptic curve that $C$ maps to, namely
\[
E_{378.a1} : z^2 = w^3 + 571293w - 68154210;
\]
we then find an explicit covering map whose $w$-coordinate is given by
\[
w(x)=\frac{9 (816 + 816 x + 1196 x^2 + 1196 x^3 + 893 x^4 + 590 x^5 + 224 x^6 + 32 x^7)}{(2 + x)^2 (8 + 3 x^2 + 2 x^3)}.
\]
Since the pullbacks to $C$ of the canonical differentials of $E_{54.a2}, E_{378.a1}$ are linearly independent, this allows us to conclude that $J$ is $(7,7)$-isogenous to the product $E_{54.a2} \times E_{378.a1}$. From this one easily concludes that $\operatorname{End}_{\mathbb{Q}}(J)=\operatorname{End}_{\overline{\mathbb{Q}}}(J)$ is an order of index 7 in $\mathbb{Z} \times \mathbb{Z}$: indeed we have proved that $\operatorname{End}_{\mathbb{Q}}(J)$ \textit{contains} an order of index 7 in $\mathbb{Z} \times \mathbb{Z}$, and on the other hand we cannot have $\operatorname{End}_{\mathbb{Q}}(J)=\mathbb{Z}^2$, because otherwise $J$ would be \textit{isomorphic} to the product of two elliptic curves, which is never the case for a Jacobian (the conclusion also depends on the fact that the curves $E_{54.a2}$ and $E_{378.a1}$ are not geometrically isogenous and do not admit potential CM -- both facts are easy to check).

\subsubsection{Quaternionic multiplication}
We study the Jacobian $J$ of the genus 2 curve
$
y^2+y = 6x^5+9x^4-x^3-3x^2;
$
this is the curve labelled 20736.l.373248.1 in the LMFDB, and the only one currently present in the database admitting quaternionic multiplication by a non-maximal order. We shall prove that $J$ does indeed admit QM over $\overline{\mathbb{Q}}$ and determine the isomorphism class of the ring $\operatorname{End}_{\overline{\mathbb{Q}}}(J)$, following the method oulined in the proof of theorem \ref{thm:CertifyQM}.

Let $x_J \in \mathcal{A}_2$ be the point in moduli space corresponding to $J$. Using the equations for the Humbert surfaces $H_{\Delta}$ for $1 \leq \Delta \leq 24$, we easily find that $x_J \in H_{12} \cap H_{24}$, and $\displaystyle x_J \not \in \bigcup_{\begin{array}{c} 1 \leq j \leq 24 \\  j \neq 12, 24 \end{array}} H_{j}$. This implies that either $J$ is $\overline{\mathbb{Q}}$-isogenous to the square of an elliptic curve or it admits quaternionic multiplication (over $\overline{\mathbb{Q}}$). In both cases, there is a quaternionic ring $R$ that acts optimally on $J_{\overline{\mathbb{Q}}}$; by theorem \ref{thm_RMImpliesEverything}, the discriminant matrix associated with this ring can be taken of the form
$
M_R = \begin{pmatrix}
12 & n \\ n & 24
\end{pmatrix},
$
where $n$ is a integer that we wish to determine, and which we can assume to be non-negative.
Combining the fact that $M_R$ is positive-definite with the fact that $\det M_R = 4\operatorname{disc}(R)$ is a multiple of 4 (see theorem \ref{thm_RMImpliesEverything}), we obtain that $n$ is even and does not exceed $16 = \lfloor \sqrt{12 \cdot 24} \rfloor$. Now notice that the quadratic form corresponding to $M_R$ primitively represents $36-2n =M_R(1,-1)$, so there is an optimal embedding of $\mathcal{O}_{36-2n}$ in $R$. On the other hand, we have already checked that $x_J$ does not belong to $H_{\Delta}$ for $0<\Delta<24$, $\Delta \neq 12$: this already implies $n \in \{0,2,4,6,12\}$. Next we check that $x_J$ does not belong to $H_{28} \cup H_{36} \cup H_{40} \cup H_{48}$, and we remark that the quadratic forms corresponding to 
\[
\begin{pmatrix}
12 & 0 \\ 0 & 24
\end{pmatrix},\begin{pmatrix} 12 & 2 \\ 2 & 24
\end{pmatrix}, \begin{pmatrix}12 & 4 \\ 4 & 24
\end{pmatrix},\begin{pmatrix} 12 & 6 \\ 6 & 24
\end{pmatrix}
\]
primitively represent $36$, $40$, $28$, and $48$ respectively (indeed they take these values on the vectors $(1,1)$, $(1,1)$, $(1,-1)$ and $(1,1)$ respectively). By the same argument as above, this implies $n \not \in \{0,2,4,6\}$, hence $n=12$. This is consistent with the LMFDB data: the endomorphism ring of $\operatorname{J}_{\overline{\mathbb{Q}}}$ was numerically computed to be an index-6 order of the unique quaternion algebra over $\mathbb{Q}$ of discriminant 6, and the endomorphism ring $R$ we just found has discriminant
$
\frac{1}{4} \operatorname{det} \begin{pmatrix}
12 & 12 \\ 12 & 24
\end{pmatrix} = 36,
$
which implies that $R$ is an order in the quaternion algebra ramified precisely at $2$ and $3$ (i.e. the quaternion algebra $Q_6$ of discriminant 6), and that its index in a maximal order is equal to $\frac{\operatorname{disc} R}{\operatorname{disc} Q_6}=6$.

\bigskip

\textbf{Acknowledgments.} I thank the authors of \cite{MR3540958}, and especially Andrew Sutherland and Jeroen Sijsling, for providing the data used for the tests of section \ref{sect:LMFDB} and for answering many questions. I also thank David Gruenewald for computing the equation of $H_{52}$. Finally, I am grateful to Christophe Ritzenthaler for many useful discussions, and to the anonymous referee for their careful reading of the manuscript and the many valuable comments.

\bibliographystyle{alpha}
\bibliography{Biblio}

\newcommand{\etalchar}[1]{$^{#1}$}
\begin{thebibliography}{{Lom}16b}

\bibitem[Ach05]{MR2181871}
Jeffrey~D. Achter.
\newblock Detecting complex multiplication.
\newblock In {\em Computational aspects of algebraic curves}, volume~13 of {\em
  Lecture Notes Ser. Comput.}, pages 38--50. World Sci. Publ., Hackensack, NJ,
  2005.

\bibitem[Ach09]{MR2496739}
J.~D. Achter.
\newblock Split reductions of simple abelian varieties.
\newblock {\em Math. Res. Lett.}, 16(2):199--213, 2009.

\bibitem[Ach12]{MR2914900}
J.~D. Achter.
\newblock Explicit bounds for split reductions of simple abelian varieties.
\newblock {\em J. Th\'eor. Nombres Bordeaux}, 24(1):41--55, 2012.

\bibitem[BLR90]{MR1045822}
S.~Bosch, W.~L{\"u}tkebohmert, and M.~Raynaud.
\newblock {\em N\'eron models}, volume~21 of {\em Ergebnisse der Mathematik und
  ihrer Grenzgebiete (3) [Results in Mathematics and Related Areas (3)]}.
\newblock Springer-Verlag, Berlin, 1990.

\bibitem[BSS{\etalchar{+}}16]{MR3540958}
Andrew~R. Booker, Jeroen Sijsling, Andrew~V. Sutherland, John Voight, and Dan
  Yasaki.
\newblock A database of genus-2 curves over the rational numbers.
\newblock {\em LMS J. Comput. Math.}, 19(suppl. A):235--254, 2016.

\bibitem[CF96]{MR1406090}
J.~W.~S. Cassels and E.~V. Flynn.
\newblock {\em Prolegomena to a middlebrow arithmetic of curves of genus
  {$2$}}, volume 230 of {\em London Mathematical Society Lecture Note Series}.
\newblock Cambridge University Press, Cambridge, 1996.

\bibitem[Chi92]{MR1156568}
W.~C. Chi.
\newblock {$l$}-adic and {$\lambda$}-adic representations associated to abelian
  varieties defined over number fields.
\newblock {\em Amer. J. Math.}, 114(2):315--353, 1992.

\bibitem[CL07]{MR2367320}
J.~E. Cremona and M.~P. Lingham.
\newblock Finding all elliptic curves with good reduction outside a given set
  of primes.
\newblock {\em Experiment. Math.}, 16(3):303--312, 2007.

\bibitem[CMSV17]{2017arXiv170509248C}
E.~{Costa}, N.~{Mascot}, J.~{Sijsling}, and J.~{Voight}.
\newblock {Rigorous computation of the endomorphism ring of a Jacobian}.
\newblock {\em ArXiv e-prints}, May 2017.

\bibitem[EK14]{MR3298543}
N.~Elkies and A.~Kumar.
\newblock K3 surfaces and equations for {H}ilbert modular surfaces.
\newblock {\em Algebra Number Theory}, 8(10):2297--2411, 2014.

\bibitem[FKRS12]{MR2982436}
Francesc Fit{\'e}, Kiran~S. Kedlaya, V{\'{\i}}ctor Rotger, and Andrew~V.
  Sutherland.
\newblock Sato-{T}ate distributions and {G}alois endomorphism modules in genus
  2.
\newblock {\em Compos. Math.}, 148(5):1390--1442, 2012.

\bibitem[Fly90]{MR1041476}
E.~V. Flynn.
\newblock The {J}acobian and formal group of a curve of genus {$2$} over an
  arbitrary ground field.
\newblock {\em Math. Proc. Cambridge Philos. Soc.}, 107(3):425--441, 1990.

\bibitem[Gon98]{MR1628150}
Josep Gonz{\'a}lez.
\newblock On the {$p$}-rank of an abelian variety and its endomorphism algebra.
\newblock {\em Publ. Mat.}, 42(1):119--130, 1998.

\bibitem[Gru08]{Gruenewald}
D.~Gruenewald.
\newblock {\em Explicit Algorithms for Humbert Surfaces}.
\newblock PhD thesis, University of Sydney, dec 2008.

\bibitem[Hal11]{MR2820155}
C.~Hall.
\newblock An open-image theorem for a general class of abelian varieties.
\newblock {\em Bull. Lond. Math. Soc.}, 43(4):703--711, 2011.
\newblock With an appendix by Emmanuel Kowalski.

\bibitem[HM95]{MR1329525}
K.~Hashimoto and N.~Murabayashi.
\newblock Shimura curves as intersections of {H}umbert surfaces and defining
  equations of {QM}-curves of genus two.
\newblock {\em Tohoku Math. J. (2)}, 47(2):271--296, 1995.

\bibitem[HS14]{MR3240808}
David Harvey and Andrew~V. Sutherland.
\newblock Computing {H}asse-{W}itt matrices of hyperelliptic curves in average
  polynomial time.
\newblock {\em LMS J. Comput. Math.}, 17(suppl. A):257--273, 2014.

\bibitem[Igu60]{MR0114819}
Jun-ichi Igusa.
\newblock Arithmetic variety of moduli for genus two.
\newblock {\em Ann. of Math. (2)}, 72:612--649, 1960.

\bibitem[Kan94]{MR1285957}
Ernst Kani.
\newblock Elliptic curves on abelian surfaces.
\newblock {\em Manuscripta Math.}, 84(2):199--223, 1994.

\bibitem[Kir69]{MR0258855}
David Kirby.
\newblock Integer matrices of finite order.
\newblock {\em Rend. Mat. (6)}, 2:403--408, 1969.

\bibitem[KM16]{MR3540944}
A.~Kumar and R.~E. Mukamel.
\newblock Real multiplication through explicit correspondences.
\newblock {\em LMS J. Comput. Math.}, 19(suppl. A):29--42, 2016.

\bibitem[Kum15]{MR3427148}
A.~Kumar.
\newblock Hilbert modular surfaces for square discriminants and elliptic
  subfields of genus 2 function fields.
\newblock {\em Res. Math. Sci.}, 2:Art. 24, 46, 2015.

\bibitem[{LMF}16]{lmfdb}
The {LMFDB Collaboration}.
\newblock The {L-functions} and modular forms database.
\newblock \url{http://www.lmfdb.org}, 2016.
\newblock [Online; accessed 07 December 2016].

\bibitem[{Lom}16a]{Surfaces}
D.~{Lombardo}.
\newblock {Explicit surjectivity for Galois representations attached to abelian
  surfaces and $\operatorname{GL}_2$-varieties}.
\newblock {\em Journal of Algebra}, 460C:26--59, 2016.

\bibitem[{Lom}16b]{CM}
D.~{Lombardo}.
\newblock {Galois representations attached to abelian varieties of CM type}.
\newblock {\em \textit{Bulletin de la Soci\'et\'e Math\'ematique de France} (to
  appear)}, 2016.

\bibitem[LP92]{MR1150604}
M.~Larsen and R.~Pink.
\newblock On {$\ell$}-independence of algebraic monodromy groups in compatible
  systems of representations.
\newblock {\em Invent. Math.}, 107(3):603--636, 1992.

\bibitem[MU01]{MR1807666}
Naoki Murabayashi and Atsuki Umegaki.
\newblock Determination of all {${\bf Q}$}-rational {CM}-points in the moduli
  space of principally polarized abelian surfaces.
\newblock {\em J. Algebra}, 235(1):267--274, 2001.

\bibitem[R{\'e}m17]{GaelUniformes}
G.~R{\'e}mond.
\newblock Conjectures uniformes sur les vari\'et\'es ab\'eliennes.
\newblock {\em The Quarterly Journal of Mathematics}, pages 1--28, 2017.

\bibitem[Riv08]{MR2401624}
Igor Rivin.
\newblock Walks on groups, counting reducible matrices, polynomials, and
  surface and free group automorphisms.
\newblock {\em Duke Math. J.}, 142(2):353--379, 2008.

\bibitem[Rot02]{RotgerThesis}
V.~Rotger.
\newblock {\em Abelian Varieties with Quaternionic Multiplication and their
  Moduli}.
\newblock PhD thesis, Universitat de Barcelona, 2002.

\bibitem[Run99]{MR1707758}
Bernhard Runge.
\newblock Endomorphism rings of abelian surfaces and projective models of their
  moduli spaces.
\newblock {\em Tohoku Math. J. (2)}, 51(3):283--303, 1999.

\bibitem[{Saw}15]{2015arXiv150604784S}
W.~F. {Sawin}.
\newblock {Ordinary Primes for Abelian Surfaces}.
\newblock {\em ArXiv e-prints}, June 2015.

\bibitem[Sil92]{MR1154704}
A.~Silverberg.
\newblock Fields of definition for homomorphisms of abelian varieties.
\newblock {\em J. Pure Appl. Algebra}, 77(3):253--262, 1992.

\bibitem[Sto95]{MR1363577}
Michael Stoll.
\newblock Two simple {$2$}-dimensional abelian varieties defined over {${\bf
  Q}$} with {M}ordell-{W}eil group of rank at least {$19$}.
\newblock {\em C. R. Acad. Sci. Paris S\'er. I Math.}, 321(10):1341--1345,
  1995.

\bibitem[vW99a]{MR1609658}
Paul van Wamelen.
\newblock Examples of genus two {CM} curves defined over the rationals.
\newblock {\em Math. Comp.}, 68(225):307--320, 1999.

\bibitem[vW99b]{MR1648415}
Paul van Wamelen.
\newblock Proving that a genus {$2$} curve has complex multiplication.
\newblock {\em Math. Comp.}, 68(228):1663--1677, 1999.

\bibitem[WM71]{MR0314847}
W.~C. Waterhouse and J.~S. Milne.
\newblock Abelian varieties over finite fields.
\newblock In {\em 1969 {N}umber {T}heory {I}nstitute ({P}roc. {S}ympos. {P}ure
  {M}ath., {V}ol. {XX}, {S}tate {U}niv. {N}ew {Y}ork, {S}tony {B}rook,
  {N}.{Y}., 1969)}, pages 53--64. Amer. Math. Soc., Providence, R.I., 1971.

\bibitem[Zar00]{MR1748293}
Yu.~G. Zarhin.
\newblock Hyperelliptic {J}acobians without complex multiplication.
\newblock {\em Math. Res. Lett.}, 7(1):123--132, 2000.

\bibitem[Zar17]{MR3660686}
Yu.~G. Zarhin.
\newblock Endomorphism rings of reductions of elliptic curves and {A}belian
  varieties.
\newblock {\em Algebra i Analiz}, 29(1):110--144, 2017.

\end{thebibliography}

\end{document}